\documentclass[11pt]{amsart}
\usepackage{amsmath}
\usepackage{amssymb}
\usepackage{mathrsfs}
\usepackage{latexsym}
\usepackage{color}
\usepackage{mathtools}
\usepackage{leftidx}
\usepackage{bm}
\usepackage{url}
\usepackage{diagrams}
\newfont{\msam}{msam10}

\newtheorem{theorem}[]{Theorem}
\newtheorem{proposition}[]{Proposition}
\newtheorem{corollary}[]{Corollary}
\newtheorem{lemma}[]{Lemma}

\theoremstyle{definition}
\newtheorem{definition}[]{Definition}
\newtheorem{defn}[theorem]{Definition}
\newtheorem{remark}[]{Remark}

\newtheorem{conj}[]{Conjecture}

\let\nc\newcommand

\def\bthm{\begin{theorem}}
\def\ethm{\end{theorem}}
\def\blemma{\begin{lemma}}
\def\elemma{\end{lemma}}
\def\bproof{\begin{proof}}
\def\eproof{\end{proof}}
\def\bprop{\begin{proposition}}
\def\eprop{\end{proposition}}
\def\bcor{\begin{corollary}}
\def\ecor{\end{corollary}}
\def\bconj{\begin{conj}}
\def\econj{\end{conj}}
\nc{\la}{\label}
\def\Q{\mathbb{Q}}

\def\L {\boldsymbol{L}}

\def\Com{\mathtt{Com}}

\def\DGL{\mathtt{DGLA}}
\def\DGC{\mathtt{DGC}}
\def\cDGC{\mathtt{DGCC}}

\def\Tw{\mathtt{Tw}}

\def\DGA{\mathtt{DGA}}

\def\cDGA{\mathtt{DGCA}}

\def\DGMod{\mathtt{DG\,Mod}}
\def\DGCoMod{\mathtt{DG\,CoMod}}
\def\DGBimod{\mathtt{DG\,Bimod}}
\def\DGBicomod{\mathtt{DG\,Bicomod}}
\def\D{\mathcal{D}}
\def\C{\mathcal{C}}

\def\U{\mathcal{U}}

\def\Ho{{\mathtt{Ho}}}
\def\mfa{\mathfrak{a}}
\nc{\ocolim}{{\rm ocolim}}
\nc{\Ob}{{\rm Ob}}
\nc{\Hom}{{\rm{Hom}}}
\nc{\RHom}{{\rm{RHom}}}
\nc{\Homcont}{{\mathcal{H}om}}
\nc{\HOM}{\underline{\rm{Hom}}}
\nc{\DER}{\underline{\rm{Der}}}
\nc{\END}{\underline{\rm{End}}}
\nc{\bSym}{\mathbf{Sym}}
\nc{\Ext}{{\rm{Ext}}}
\nc{\Rep}{{\rm{Rep}}}
\nc{\DRep}{{\rm{DRep}}}
\nc{\ODRep}{{\mathcal O}{\rm{DRep}}}
\nc{\NCRep}{\widetilde{\rm{Rep}}}
\nc{\RAct}{{\rm{RAct}}}
\nc{\bs}{\backslash}
\nc{\ob}{{\tt{Obs}}}
\nc{\CE}{\mathcal{C}}
\nc{\TP}{{T\!P}}
\nc{\un}{\underline{n}}
\nc{\um}{\underline{m}}
\nc{\rn}{\langle n \rangle}
\nc{\nn}{{{\natural} {\natural}}}
\nc{\n}{{{\natural}}}
\nc{\A}{\mathbb A}
\nc{\B}{{\mathrm{B}}}
\nc{\Ba}{\overline{\mathrm{B}}}
\nc{\bC}{\overline{C}}
\nc{\bOmega}{\boldsymbol{\Omega}}
\nc{\bB}{\boldsymbol{B}}
\nc{\EXT}{\underline{\rm{Ext}}}
\nc{\TOR}{\underline{\rm{Tor}}}
\def\H{\mathrm H}

\def\HC{\mathrm{HC}}
\def\HR{\mathrm{HR}}
\def\rHC{\overline{\mathrm{HC}}}
\def\rHH{\overline{\mathrm{HH}}}

\nc{\End}{{\rm{End}}}
\nc{\GL}{{\rm{GL}}}
\nc{\gl}{{\mathfrak{gl}}}
\nc{\rgl}{\overline{{\mathfrak{gl}}}}
\nc{\g}{{\mathfrak{g}}}
\nc{\h}{{\mathfrak{h}}}
\nc{\PGL}{{\rm{PGL}}}
\nc{\SL}{{\rm{SL}}}
\nc{\sll}{\mathfrak{sl}}
\nc{\cn}{ \mbox{\rm c\^{o}ne} }
\nc{\PSL}{{\rm{PSL}}}
\nc{\ad}{{\rm{ad}}}
\nc{\Ad}{{\rm{Ad}}}
\nc{\dlim}{\varinjlim}
\nc{\plim}{\varprojlim}
\nc{\colim}{{{\rm colim}}}

\newcommand{\HH}{{\rm{HH}}}

\newcommand{\Sym}{\mathrm{Sym}}

\newcommand{\id}{{\rm{Id}}}
\newcommand{\Der}{{\rm{Der}}}

\newcommand{\Tr}{{\rm{Tr}}}

\newcommand{\Ker}{{\rm{Ker}}}

\def\cb{\boldsymbol{\Omega}}

\def\bs{\backslash}

\def\Gr{\mathtt{Gr}}

\def\ffgr{\mathfrak{G}}

\def\LL{\mathcal{L}}

\newcommand{\rar}{\xrightarrow{}}

\nc{\env}{\mathrm{End}(V)}
\nc{\FT}{\mathcal{C}}

\numberwithin{equation}{section}
\numberwithin{theorem}{section}
\numberwithin{lemma}{section}
\numberwithin{proposition}{section}
\numberwithin{definition}{section}
\numberwithin{corollary}{section}
\numberwithin{example}{section}
\numberwithin{remark}{section}

\newcommand{\rH}{\overline{\mathrm{H}}}

\def\bdf{\begin{defn}}
\def\edf{\end{defn}}

\def\brm{\begin{remark}}
\def\erm{\end{remark}}

\theoremstyle{definition}
\def\bdf{\begin{definition}}
\def\edf{\end{definition}}

\def\arbreBA{\vcenter{\xymatrix@R=2pt@C=2pt{
&&&&\\
&&&*{}\ar@{-}[ul] & \\
&&*{}\ar@{-}[uurr] \ar@{-}[uull] \ar@{-}[d]     &&\\
&&&&
}}}

\def\arbreAB{\vcenter{\xymatrix@R=2pt@C=2pt{
&&&&\\
&*{}\ar@{-}[ur] &&& \\
&&*{}\ar@{-}[uurr] \ar@{-}[uull] \ar@{-}[d]     &&\\
&&&&
}}}

\def\arbreABC{\vcenter{\xymatrix@R=1pt@C=1pt{
&&&&&&\\
&*{}\ar@{-}[ur] &&&&& \\
&&*{}\ar@{-}[uurr] &&&&\\
&&&*{}\ar@{-}[uuurrr] \ar@{-}[uuulll] \ar@{-}[d] &&&\\
&&&&&&
}}}

\def\arbreBAC{\vcenter{\xymatrix@R=1pt@C=1pt{
&&&&&&\\
&&&*{}\ar@{-}[ul] &&& \\
&&*{}\ar@{-}[uurr] &&&&\\
&&&*{}\ar@{-}[uuurrr] \ar@{-}[uuulll] \ar@{-}[d] &&&\\
&&&&&&
}}}

\def\arbreACB{\vcenter{\xymatrix@R=1pt@C=1pt{
&&&&&&\\
&*{}\ar@{-}[ur] &&&&& \\
&&&&*{}\ar@{-}[uull] &&\\
&&&*{}\ar@{-}[uuurrr] \ar@{-}[uuulll] \ar@{-}[d] &&&\\
&&&&&&
}}}

\def\arbreBCA{\vcenter{\xymatrix@R=1pt@C=1pt{
&&&&&&\\
&&&&&*{}\ar@{-}[ul] & \\
&&*{}\ar@{-}[uurr] &&&&\\
&&&*{}\ar@{-}[uuurrr] \ar@{-}[uuulll] \ar@{-}[d] &&&\\
&&&&&&
}}}

\def\arbreCAB{\vcenter{\xymatrix@R=1pt@C=1pt{
&&&&&&\\
&&&*{}\ar@{-}[ur] &&& \\
&&&&*{}\ar@{-}[uull] &&\\
&&&*{}\ar@{-}[uuurrr] \ar@{-}[uuulll] \ar@{-}[d] &&&\\
&&&&&&
}}}

\def\arbreCBA{\vcenter{\xymatrix@R=1pt@C=1pt{
&&&&&&\\
&&&&&*{}\ar@{-}[ul] & \\
&&&&*{}\ar@{-}[uull] &&\\
&&&*{}\ar@{-}[uuurrr] \ar@{-}[uuulll] \ar@{-}[d] &&&\\
&&&&&&
}}}

\def\arbreACA{\vcenter{\xymatrix@R=1pt@C=1pt{
&&&&&&\\
&*{}\ar@{-}[ur] &&&&*{}\ar@{-}[ul] & \\
&&&&&&\\
&&&*{}\ar@{-}[uuurrr] \ar@{-}[uuulll] \ar@{-}[d] &&&\\
&&&&&&
}}}

\begin{document}

\title{Hodge decomposition of string topology}
\author{Yuri Berest}
\address{Department of Mathematics,
Cornell University, Ithaca, NY 14853-4201, USA}
\email{berest@math.cornell.edu}
\author{Ajay C. Ramadoss}
\address{Department of Mathematics,
Indiana University,
Bloomington, IN 47405, USA}
\email{ajcramad@indiana.edu}
\author{Yining Zhang}
\address{Department of Mathematics,
University of Colorado Boulder,
Boulder, CO 80309, USA}
\email{yining.zhang@colorado.edu}

\begin{abstract}
Let $X$ be a simply connected closed oriented manifold of rationally elliptic homotopy type.
We prove that the string topology bracket on the $S^1$-equivariant homology
$ \rH_\ast^{S^1}(\LL X,\Q) $ of the free loop space of $X$ preserves the Hodge decomposition
of $ \rH_\ast^{S^1}(\LL X,\Q) $, making it a bigraded Lie algebra. We deduce this result from a general
theorem on derived Poisson structures on the
universal enveloping algebras of homologically nilpotent finite-dimensional DG Lie algebras.
Our theorem settles a conjecture of \cite{BRZ}.
\end{abstract}

\maketitle

\section{Introduction}
Let $X$ be a simply connected closed oriented manifold, and let $\LL X :=\mathrm{Map}(S^1,X)$ denote the free loop space over $X$. Chas and Sullivan \cite{ChS} showed that $ \rH_\ast^{S^1}(\LL X,\Q)$, the (reduced) rational $ S^1$-equivariant homology of $\LL X $ with respect to the natural circle action, carries a graded Lie algebra structure with the so-called {\it string topology bracket}
\begin{equation} \la{stringtop} \{\, \mbox{--}\,,\,\mbox{--}\, \}:\ \rH_\ast^{S^1}(\LL X, \Q) \times \rH_\ast^{S^1}(\LL X, \Q) \rar \rH_\ast^{S^1}(\LL X, \Q)\ . \end{equation}
This bracket is intrinsically related to the geometry of $ \LL X $ and has many interesting properties, which have been studied extensively in recent years (see, e.g., \cite{FeT,FTV,TZ1,TZ2}).

In this paper we show that \eqref{stringtop} is compatible with Frobenius (power) operations on $\rH_\ast^{S^1}(\LL X, \Q)$: i.e., the string topology bracket respects a Hodge-type decomposition:
\begin{equation}
\la{hodgedecompos} \rH_\ast^{S^1}(\LL X, \Q)\, = \,
\bigoplus_{p=0}^{\infty}\, \rH_\ast^{S^1,\,(p)}(\LL X, \Q)\,,
\end{equation}
where the direct summands are common eigenspaces of graded endomorphisms of
$\rH_\ast^{S^1}(\LL X, \Q)$ with eigenvalues $n^p$, $\,n\ge 0$, induced by the finite coverings of the
circle: $\,S^1 \rar S^1\,,\,\, e^{i\theta} \mapsto e^{ni\theta}$ (see \cite{BFG91}).
%
%
%
More precisely, we prove the following theorem.
\bthm 
\la{MainTheorem}
Assume that the manifold $X$ is rationally elliptic as a topological space, that is\\
$\, \dim\, \sum_{i\ge 2} \pi_{i}(X) \otimes \Q < \infty \,$.
Then
$$
\{\rH_\ast^{S^1,\,(p)}(\LL X, \Q),\, \rH_\ast^{S^1,\,(q)}(\LL X, \Q) \} \,\subseteq \, \rH_\ast^{S^1,\,(p+q-1)}(\LL X, \Q) \ ,\quad \forall\,p,q \ge 0\, ,\ p+q \ge 1\ .
$$
Thus, the Chas-Sullivan Lie algebra of $X$ is bigraded:
$$
\rH_\ast^{S^1}(\LL X,\Q) \, = \, \bigoplus_{n\ge 0}\,
\bigoplus_{p \ge -1}\ \rH_n^{S^1,\,(p+1)}(\LL X, \Q)\,,
$$
where the first grading is given by homological degree and the second by Hodge
degree $($shifted by one$)$.
\ethm

The result of Theorem \ref{MainTheorem} was conjectured in our earlier paper \cite{BRZ}, where we proved that the bracket \eqref{stringtop} preserves a filtration on the vector space $\rH^{S^1}_\ast(\LL X, \Q)$ associated naturally with the direct sum decomposition \eqref{hodgedecompos}. Thus we strengthen the main result of \cite{BRZ}, albeit under the additional assumption that $X$ is rationally elliptic\footnote{Recall that,
in rational homotopy theory, there is a well-known dichotomy dividing all simply connected spaces with finite rational homology into two classes: elliptic and hyperbolic. Although `generic' spaces are
known to be rationally hyperbolic, many important spaces occurring `in nature' are rationally elliptic: these include, for example, the spheres $ S^n $ ($n\ge 2$), the complex projective spaces
$ \mathbb{C} \mathbb{P}^r $ $(r \ge 1)$, all compact connected Lie groups $G$ and
their homogeneous spaces $G/K$ with $K $ compact connected. Moreover, any simply connected compact manifold $X$ of dimension $ d $ is known to be rationally elliptic if
$ \pi_i(X) \otimes \Q = 0 $ for $ i > d $ (see \cite[Part IV, Sect. 32]{FHT}).}.

Theorem \ref{MainTheorem} is a geometric fact: it relates two geometrically defined structures on a smooth manifold $X$. Unfortunately, we do not know how to see this relation directly, in geometric terms, using the original definition of string topology in \cite{ChS}. Instead, we prove Theorem \ref{MainTheorem} in a somewhat roundabout way, deducing it from an abstract algebraic result on derived Poisson structures on the universal enveloping algebra $\U\mfa$ of a (DG) Lie algebra $\mfa$ (see \cite{BCER}). The main property of such a structure is that it induces naturally a  Lie bracket on the (reduced) cyclic homology of $\U\mfa$:
\begin{equation} \la{dpbracket}
\{\, \mbox{--}\,,\,\mbox{--}\,\}:\ \rHC_\ast(\U\mfa) \times \rHC_\ast(\U\mfa) \,\rar\, \rHC_\ast(\U\mfa)
\end{equation}
which is an algebraic model for the string topology bracket \eqref{stringtop}. On the other hand, for any DG Lie algebra $\mfa$, the cyclic homology of $\U\mfa$ has a canonical direct sum decomposition
\begin{equation}
\la{hodgedsIntro}
\rHC_{\bullet}(\U\mfa) \,=\,\bigoplus_{p=1}^{\infty}\, \HC^{(p)}_{\bullet}(\mfa)
\end{equation}
which is called the {\it Lie-Hodge decomposition} (see \cite{BFPRW, BRZ} and also Section \ref{secHodgeDecomp} below). This raises a natural question about compatibility of the two structures: 
namely,

\begin{center}
Does the derived Poisson bracket \eqref{dpbracket} preserve \eqref{hodgedsIntro} ?
\end{center}

As shown in \cite{BRZ}, the answer to this question is, in general, negative. It is therefore
necessary to impose certain restrictions on the Lie algebra $ \mfa $ and the derived Poisson structure on $\U\mfa$. In this paper, we consider a special class of derived Poisson structures on $\U\mfa$ that arise from a cyclic pairing on a cocommutative DG coalgebra $C$ Koszul dual to the DG Lie algebra $\mfa$. Then, under natural finiteness assumptions on $ \mfa $ and $ C $, we prove that \eqref{dpbracket} does preserve \eqref{hodgedsIntro}. 
To state our main result in precise terms, we recall that every DG Lie algebra $\mfa$ has a (unique) minimal model, which is given by an $L_\infty$-algebra structure on the homology $\H_\ast(\mfa)$ of $ \mfa $ together with a canonical $L_\infty$-quasi-isomorphism $\mfa \stackrel{\sim}{\rar} \H_\ast(\mfa)$. We denote this minimal $L_\infty$-model simply by $\H_\ast(\mfa)$.
\bthm 
\la{genmain}
Let $\mfa\,\in\,\DGL^{+}_k$  be a non-negatively graded DG Lie algebra defined over a field $k$ of characteristic $0$. Assume that\\
$(1)$ $\dim_k \H_\ast(\mfa)<\infty$  and $\H_\ast(\mfa)$ is nilpotent as an $L_\infty$ algebra.\\
$(2)$ $\mfa$ has a Koszul dual cocommutative coalgebra $C$ of finite total dimension, i.e. $\dim_k C <\infty$.\\
Then the derived Poisson bracket \eqref{dpbracket} associated to a(ny) nondegenerate cyclic pairing on $C$ preserves the Hodge decomposition \eqref{hodgedsIntro}: i.e.,
$$ \{\HC^{(p)}_\ast(\mfa)\,,\, \HC^{(q)}_\ast(\mfa)\} \,\subseteq\,\HC^{(p+q-2)}_\ast(\mfa) \ ,
\quad \forall\,p,q \ge 0 \ .
$$
\ethm
Note that Theorem~\ref{genmain} applies, in particular, to an ordinary finite-dimensional nilpotent Lie algebra $ \mfa $, with the derived Poisson bracket on $ \rHC_\ast(\U\mfa) $ coming from the
natural pairing on the Chevalley-Eilenberg chain complex $ \,\C_\ast(\mfa;k) = \wedge^{\ast} \mfa \,$ (see \cite[Sect. 6]{CEEY}).

Now, Theorem \ref{MainTheorem} is a consequence of Theorem \ref{genmain} modulo known results in the literature. First, if $X$ is a simply connected manifold, we take $\mfa = \mfa_X $ to be its Quillen model \cite{Q2}. Then, by a classical theorem of Goodwillie \cite{Go} and Jones \cite{J}, there is a natural isomorphism $\rHC_\ast(\U\mfa)\,\cong\,\rH_\ast^{S^1}(\LL X, \Q)$; moreover, as shown in \cite[Theorem 1.2]{BRZ}, this isomorphism identifies $\,\HC^{(p)}(\mfa)\,\cong\,\rH_\ast^{S^1,(p-1)}(\LL X, \Q)$ for all $p \geq 1$. Thus, the geometric Hodge decomposition \eqref{hodgedecompos} for a simply connected space $X$ coincides (up to a shift in degree) with the Lie-Hodge decomposition \eqref{hodgedsIntro} for the Lie model $ \mfa $ of $X$. Next, for any compact manifold $X$,  Lambrechts and Stanley \cite{LS} constructed a finite-dimensional commutative algebra model $A_X$, whose (linear) dual coalgebra $ C := \Hom(A_X,  \Q) $ is Koszul dual to the Quillen model $\mfa$. This coalgebra $C$ comes equipped with a nondegenerate cyclic pairing (Poincar\'{e} duality), and --- as observed in \cite{BCER} (see also Lemma \ref{StringPoiss} below) ---  the associated derived Poisson bracket on $\rHC_\bullet(\U\mfa)$ agrees with the Chas-Sullivan  bracket on $\rH_\ast^{S^1}(\LL X, \Q)$. To apply Theorem \ref{genmain} it remains note that
for $X$ rationally elliptic, the minimal $L_\infty$-model of $\mfa_X$ is finite-dimensional and nilpotent: i.e., $\H_\ast(\mfa)$ satisfies condition $(1)$ of Theorem \ref{genmain}.

Next, we briefly outline our proof of Theorem \ref{genmain}. As a first step, we replace the cyclic homology $\rHC_\ast(\U\mfa)$ of the algebra $ \U \mfa $ by its Hochschild cohomology $\HH^{\ast}(\U\mfa,\U\mfa)$ and, following the idea of  \cite{CEEY}, express the derived Poisson bracket on $\rHC_\ast(\U\mfa)$ in terms of the canonical cup product on $\HH^{\ast}(\U\mfa,\U\mfa)$ (see Proposition \ref{poissoncup}). We show that the Lie-Hodge decomposition  of $\rHC_\ast(\U\mfa)$  naturally extends to a direct sum decomposition of $\HH^{\ast}(\U\mfa,\U\mfa)$ which we also refer to as a Lie-Hodge decomposition (see Theorem \ref{dualHodgeCoh}). Using the results of \cite{BRZ} and \cite{CEEY}, we then reduce the proof of Theorem \ref{genmain} to proving that the cup product on $\HH^{\ast}(\U\mfa,\U\mfa)$ preserves its Lie-Hodge decomposition. Writing  $\mathcal{A}$ for the Chevalley-Eilenberg cochain complex of the minimal $L_\infty$-model $ \H_*(\mfa) $ of the DG Lie algebra $\mfa$, we observe that there is a (Hodge degree preserving) algebra isomorphism $\,\HH^{\ast}(\U\mfa,\U\mfa)\,\cong\,\HH^{\ast}_{\oplus}(\mathcal{A},\mathcal{A})\,$, where $ \HH^{\ast}_{\oplus}(\mathcal{A},\mathcal{A}) $ stands for the Hochschild cochain complex of $ {\mathcal A} $ constructed using direct sums (in place of infinite direct products). Thus, our problem reduces to showing that the cup product preserves the Lie-Hodge decomposition of $\,\HH^{\ast}_{\oplus}(\mathcal{A},\mathcal{A})$. To prove this we recall that, by definition, $\,\mathcal{A}$ is a symmetric algebra equipped with a differential which encodes the $L_\infty$ structure on $\,\H_*(\mfa) $. Hence, by the Hochschild-Kostant-Rosenberg Theorem, there is a natural linear map $\,\mathrm{I}_{\mathrm{HKR}}:\,\H^{\ast}[\mathcal{V}] \rar \HH^{\ast}_{\oplus}(\mathcal{A},\mathcal{A})\,$, where $\H^{\ast}[\mathcal{V}]$  denotes the cohomology ring of the algebra $\mathcal{V} = \mathcal{V}(\mathcal{A})$  of poly-derivations of $ \mathcal{A} $. It is easy to show that the map $ \mathrm{I}_{\mathrm{HKR}} $ preserves Hodge grading; however, in general, it is {\it not} a homomorphism of graded algebras. In fact, Kontsevich's (cohomological) version of Duflo's classical theorem, which applies to $\mathcal{A}$ because $\H_*(\mfa) $ is finite-dimensional, says that  to obtain an algebra isomorphism between $\H^{\ast}[\mathcal{V}]$ and $\HH^{\ast}_{\oplus}(\mathcal{A},\mathcal{A})$ one needs to `correct' the HKR map by contracting it with a certain canonical cohomology class called the {\it Todd genus}. In general, the contraction by the Todd genus does not preserve the Hodge grading on $\mathcal{V}$; however, when $\H_*(\mfa) $ is a nilpotent $L_\infty$-algebra, we show (see Proposition \ref{DufloRE}) that the Todd genus is actually trivial\footnote{In the case when
$ \mfa $ is an ordinary finite-dimensional nilpotent Lie algebra, this observation goes back essentially to Duflo's original paper \cite{Du} (see also \cite{PT}).}. Thus, it turns out that,  under our assumptions on $ \mfa $, the HKR map is an isomorphism of graded algebras. Now, since the product on $\mathcal{V}$ obviously preserves the Hodge grading, this completes the proof of Theorem \ref{genmain}.

We would like to conclude this introduction by mentioning a well-known analogy between 
rational homotopy theory and local commutative algebra (see \cite{AH}). The rationally
elliptic spaces correspond in commutative algebra to local complete intersection rings, and
the $S^1$-equivariant homology of free loop spaces correspond to (relative) cyclic homology of
local rings. It seems natural to ask if there is a result parallel to our Theorem~\ref{MainTheorem}
in commutative algebra. Although there are well-defined Lie algebra models associated to local rings, 
our ``abstract'' Theorem~\ref{genmain} cannot be applied directly to such models, since they are no 
longer defined over characteristic zero fields.

The paper is organized as follows. In Section \ref{secHodgeDecomp}, we review the Loday-Goodwillie
 (simplicial) approach to Hodge decompositions and show how the Lie-Hodge decomposition of $\rHC_\ast(\U\mfa)$ constructed in \cite{BFPRW} arises from this approach. We also construct the Lie-Hodge decompositon on the Hochschild {\it cohomology} $\HH^{\ast}(\U\mfa,\U\mfa)$ and prove some technical results needed for our main theorem. In Section \ref{seccyclicpoiss}, we recall definitions and review some known results on derived Poisson structures. The most important for us result (proven in \cite{CEEY}) relates the cyclic derived Poisson structure on $ \U\mfa $ to the cup product and Gerstenhaber bracket on the Hochschild cohomology of $ \U \mfa $ via the Van den Bergh duality. We also show that the Van den Bergh duality is compatible with Lie-Hodge decompositions (Lemma \ref{hodgeduality}). Finally, Section \ref{SecDerivedPoissonHodge} contains the proof of Theorem \ref{genmain} as well as its application to string topology.

\subsection*{Acknowledgements}
The first author is grateful to the Department of Mathematics of Indiana
University (Bloomington) and Forschungsinstitut f\"{u}r Mathematik (ETH,
Z\"{u}rich) for their hospitality and financial support during his stay in Fall 2019. The research of the first author was partially supported by
2019 Simons Fellowship and NSF grant DMS 1702372. Research of the second author was partially supported by NSF grant DMS 1702323.

\section{Hodge decompositions}
\la{secHodgeDecomp}
It is well known that the cyclic homology of any commutative (DG) algebra $A$ has a natural decomposition
\begin{equation}
\la{hcyccomm}
\HC_\ast(A)\,=\, \bigoplus_{p=0}^{\infty} \HC^{(p)}_\ast(A)\ .
\end{equation}
which is usually called the {\it Hodge $($or $\lambda$-$)$decomposition} of $ \HC_\ast(A) $. J.-L. Loday \cite{L89} gave an elegant explanation of this phenomenon in terms of the classical bar construction
$\mathrm{C}_\ast(A)$ of the algebra $A$. Recall that for any associative algebra, $\mathrm{C}_\ast(A)$ is a cyclic module: that is, a contravariant functor ${\Delta}C^{\mathrm{op}} \rar \Com_k$ defined on the Connes cyclic category ${\Delta}C$. The category $\Delta C^{\mathrm{op}}$ naturally embeds into the category of finite sets, $\mathbf{Fin}$, and the theorem of Loday ({\it cf}. \cite[6.4.4]{L}) asserts that the cyclic homology of a cyclic module $E_\ast\,:\,\Delta C^{\mathrm{op}} \rar \Com_k$ admits a direct sum decomposition \eqref{hcyccomm} whenever the functor $E_{\ast}$ extends to $\mathbf{Fin}$, i.e. factors through the inclusion $\Delta C^{\mathrm{op}} \hookrightarrow \mathbf{Fin}$. Now, for the cyclic bar construction $E_\ast=\mathrm{C}_\ast(A)$, this happens exactly when $A$ is a commutative algebra.

T. Goodwillie ({\it cf}. \cite[6.4.5]{L}) observed that the Hodge decomposition of
$ \HC_\ast(E) $ exists in a more general situation: namely, when the functor
$ E_\ast:\,\Delta C^{\mathrm{op}} \rar \Com_k $ admits an extension $ E^{\Psi}_\ast:\,
\Delta \Psi^{\mathrm{op}} \rar \Com_k \,$ to the so-called
{\it epicyclic category} $\, \Delta \Psi^{\mathrm{op}} $. The category $
\Delta \Psi $ has the same objects as $ \Delta C $ but contains extra morphisms which
induce the power (Adams) operations on $ \HC_\ast(E) $; however, $ \Delta \Psi $ is strictly smaller
than $\mathbf{Fin}$. Thus, if $ E_\ast:\,\Delta C^{\mathrm{op}} \rar \Com_k$ factors through the inclusion $\Delta C^{\mathrm{op}} \hookrightarrow \Delta \Psi^{\mathrm{op}}$, then $\HC_\ast(E)$ has a natural Hodge decomposition; furthermore, this decomposition agrees with Loday's when the corresponding epicyclic module $E^{\Psi}_\ast\,:\,\Delta \Psi^{\mathrm{op}} \rar \Com_k$ factors through  $\Delta \Psi^{\mathrm{op}} \hookrightarrow \mathbf{Fin}$.

In this section we show that the cyclic bar construction $\mathrm{C}_\ast(\U\mfa)$ of the universal enveloping algebra of a (DG) Lie algebra $\mfa$ has a natural epicyclic structure. We prove that the Hodge decomposition of cyclic homology $\rHC_\ast(\U\mfa)$ arising from this epicyclic structure coincides with the Lie-Hodge decomposition constructed in \cite{BFPRW}. We will also establish some properties of Hodge decomposition of Hochschild homology which we will need for the proof of our main theorem.

\subsection{Epicyclic modules and Adams operations} \la{secepicyclic}

Let $\Delta$ denote the standard simplicial category whose objects are the finite ordered sets $[n]\,=\,\{0<1<2<\ldots<n\}$, and morphisms are the order preserving maps. It is easy to prove that $\Delta$ is generated by two families of maps: $\{d^i_n\,:\,[n-1] \rar [n]\}_{0 \leq i \leq n, n \geq 1}$ and $\{s^j_n \,:\,[n+1] \rar [n]\}_{0 \leq j \leq n, n \geq 0}$, called the (co)face and (co)degeneracy maps, respectively. These maps satisfy the standard (co)simplicial relations which are given, for example, in \cite[Appendix B.3]{L}. Connes' {\it cyclic category} $\Delta C$ is an extension of $\Delta$ that contains --- in addition to the $d^i_n$'s and $s^j_n$'s ---  the cyclic maps $\{\tau_n\,:\,[n] \rar [n]\}_{n \geq 0}$. More generally, for any integer $k \geq 1$, we can define the $k$-{\it cyclic category} $\Delta C^{(k)}$ that contains $ \Delta $ (and has the same objects as $\Delta$) with additional morphisms ${}_k\tau_n\,:\,[n] \rar [n]$ satisfying
$$
{}_k\tau_n \circ d^i_n\,=\, d^{i-1}_n{}_k \circ \tau_n\,, \qquad {}_k\tau_n \circ s^j_n\,=\, s^{j-1}_n{}_k\circ \tau_n\,, \qquad ({}_k\tau_n)^{k(n+1)}=\id_{[n]} \ .
$$
There are two natural functors relating $\Delta C^{(k)}$ to $ \Delta C^{(1)} \equiv \Delta C \,$:
\begin{equation} \la{powersd}
\mathrm{P}^k\,:\,\Delta C^{(k)} \rar \Delta C\,, \qquad \mathrm{Sd}^k\,:\,\Delta C^{(k)} \rar \Delta C\ .\end{equation}
The functor $\mathrm{P}^k$ is characterized by the property that its restriction to $\Delta$ is the identity,  while $\mathrm{P}^k({}_k\tau_n)=\tau_n$ for all  $n \geq 0$. The functor $\mathrm{Sd}^k$ --- called the {\it $k$-th edgewise subdivision functor} --- is defined  by
$$ \mathrm{Sd}^k([n]):=[k(n+1)-1]\,=\,[n] \,\sqcup \,\stackrel{k}{\ldots} \,\sqcup \,[n]\,,$$
and on morphisms:
$$  \mathrm{Sd}^k(\varphi)\,=\, \varphi\, \sqcup \,\stackrel{k}{\ldots} \,\sqcup\, \varphi\,,\text{ for } \varphi \,\in\,\mathrm{Mor}(\Delta)\,, \qquad  \mathrm{Sd}^k({}_k\tau_n)\,=\,\tau_{k(n+1)-1} \ . $$
Now, the {\it epicyclic category} $\Delta \Psi$ (see \cite{BFG}) is the extension of $\Delta C$ (i.e, $\Delta C \subset \Delta \Psi$), which --- in addition to the morphisms $\{d^i_n\}, \{s^j_n\}$ and $\{\tau_n\}$ generating $\Delta C$ --- contains a family of morphisms
$$ \pi_n^k:\,[k(n+1)-1] \rar [n]\,,\quad \forall n \geq 0\, ,\ k \geq 1\,, $$
called the {\it $($co$)$power maps}. These maps are characterized by the property that $\pi^k_\ast\,:\,\mathrm{Sd}^k \rar \mathrm{P}^k$ define natural transformations of functors $\Delta C^{(k)} \rar \Delta C \hookrightarrow \Delta \Psi$ for all $k \geq 1$, and, in addition, satisfy the relations
$$ \pi_n^1\,=\,\id_{[n]}\,,\qquad \pi_n^l \circ \pi^k_{l(n+1)-1}\,=\,\pi_n^{kl}\ .$$
If $\C$ is a category, an epicyclic object in $\C$ is, by definition, a functor $X: \Delta \Psi^{\mathrm{op}} \rar \C$; we will write $\C_{\Delta \Psi}$ for the category of such functors, with morphisms being the natural transformations. Note that giving an epicyclic object in $\C$ is equivalent to giving a cyclic object $X:\Delta C^{\mathrm{op}} \rar \C$ together with a family of morphisms in $\Delta C^{(k)}$ ($k \geq 1$):
\begin{equation} \la{powerops}
p^k_\ast(X):\,\mathrm{P}_*^k(X) \rar \mathrm{Sd}_*^k(X)\,,
\qquad p^k_n\,:\,X_n \rar X_{k(n+1)-1}\,,
\end{equation}
satisfying
$$ p^1_\ast=\id\ \,,\qquad p^k_\ast \circ p^l_\ast=p^{kl}_\ast $$
where $\mathrm{P}_*^k(X)$ and $\mathrm{Sd}_*^k(X)$ are the $k$-cyclic objects $(\Delta C^{(k)})^{\mathrm{op}} \rar \C$ defined by $\mathrm{P}_*^k(X) := X \circ \mathrm{P}^k $ and $\mathrm{Sd}_*^k(X):=X \circ \mathrm{Sd}^k$ (see \eqref{powersd}). It is a classical observation (due to A. Connes) that when $\C=\mathbf{Set}$, the geometric realization $|X|$ of any cyclic set $X\,:\,\Delta C^{\mathrm{op}} \rar \mathbf{Set}$ caries a natural $S^1$-action (see \cite[Thm. 7.1.4]{L}). In a similar way, if $X$ is an epicyclic set then in addition to the $S^1$-action, its realization $|X|$ carries power operations which induce Adams operations on the $S^1$-equivariant homology $\H^{S^1}_\ast(|X|)$ of $|X|$ (see \cite[Theorem A]{BFG}). We will look at an algebraic (chain) version of this construction.

Recall that to any associative unital (DG) $k$-algebra $A$, one can attach naturally a cyclic module
\begin{equation} \la{cyclicbar} \mathrm{C}_\ast(A)\,:\,\Delta C^{\mathrm{op}} \rar \Com_k\,, \, [n] \mapsto A^{\otimes (n+1)} \end{equation}
called the {\it cyclic bar construction} (see \cite[2.5.4]{L}). Our main observation in this section is
\bprop \la{cHopfepi}
If $A$ is a cocommutative (DG) Hopf algebra, then $\mathrm{C}_\ast(A)$ carries a natural epicyclic structure.
\eprop

To prove Proposition \ref{cHopfepi} we need to extend the functor \eqref{cyclicbar} to the epicyclic category, i.e. to construct a functor $\mathrm{C}^{\Psi}_{\ast}(A)\,:\,\Delta \Psi^{\mathrm{op}} \rar \Com_k$ such that $\mathrm{C}^{\Psi}_{\ast}(A)|_{\Delta C}=\mathrm{C}_{\ast}(A)$. This can be done directly by defining the structure maps \eqref{powerops} in an explicit way and verifying the required relations. We will give a more conceptual (categorical) construction of $\mathrm{C}^{\Psi}_{\ast}(A)$. To this end we will use a known characterization of cocommutative (DG) Hopf algebras as group objects in $\Com_k$. Let $\ffgr$ denote (the skeleton of) the category of finitely generated free groups: thus, the objects of $\ffgr$ are the free groups $\langle n\rangle=\mathbb{F}_n$, one for each cardinality $n \geq 0$, and the morphisms $\langle n \rangle \rar \langle m \rangle$ are arbitrary group homomorphisms $\mathbb{F}_n \rar \mathbb{F}_m$. The category $\ffgr$ carries a (strict) monoidal structure with product $\langle n \rangle \ast \langle m \rangle = \langle n+m\rangle$ for all $n,m \geq 0$. The category of all (discrete) groups $\Gr$ can then be described as the category $\mathbf{Set}^{\otimes}_{\ffgr}$ of strict monoidal functors $\ffgr^{\mathtt{op}} \rar \mathbf{Set}$ with values in $\mathbf{Set}$ equipped with the usual (cartesian) monoidal structure:  the equivalence $\Gr \stackrel{\sim}{\rar} \mathbf{Set}^{\otimes}_{\ffgr}$ is given by the Yoneda functor $G \mapsto \underline{G}:=\Hom_{\Gr}(\mbox{--},G)$ restricted to the subcategory $\ffgr \subset \Gr$. Now, it is known and easy to check (see, e.g., \cite{BHM} and \cite[Example 1.3]{BFG}) that the cyclic nerve $N^{\mathrm{cyc}}_\ast(G):=\{G^{n+1}\}_{n \geq 0}$ of any discrete group $G$ carries a canonical epicyclic structure with power maps \eqref{powerops} given by
$$
p^k_n\,:\,G^{n+1} \rar G^{k(n+1)} \,,\qquad (g_0,\cdots,g_n) \mapsto (g_0,\cdots, g_n; \stackrel{k}{\ldots}; g_0,\cdots,g_n)\ .
$$
Thus we have a well-defined functor $N^{\mathrm{cyc}}_\ast:\Gr \rar \mathbf{Set}_{\Delta \Psi}$.  If we identify $\Gr \,\cong\, \mathbf{Set}^{\otimes}_{\ffgr}$ via $G \mapsto \underline{G}$ as above,  then
$N^{\mathrm{cyc}}_\ast $ is simply the pull-back functor $\Psi^{\ast}:\, \mathbf{Set}^{\otimes}_{\ffgr} \rar \mathbf{Set}_{\Delta \Psi}$ for a natural map
 \begin{equation} \label{epiintog} \Psi\,:\,\Delta \Psi \rar \ffgr \ .\end{equation}
Explicitly, \eqref{epiintog} is defined  on objects by $\Psi([n])=\langle n+1 \rangle=\mathbb{F}\langle x_0,\ldots,x_n\rangle$ and on morphisms by the following  formulas:
\begin{equation*}
\left. \begin{aligned}
\Psi(d^i_n) &: \langle n \rangle \rar \langle n+1 \rangle\,,\qquad (x_0, x_1, \ldots, x_{n-1}) \mapsto \begin{cases} (x_0,\ldots, x_{i-1}, x_i x_{i+1}, \ldots, x_n)\ , &  0 \leq i < n\\
(x_nx_0, x_1,\ldots, x_{n-1})\ , & i=n \end{cases} \\
\Psi(s^j_n) &: \langle n+2\rangle \rar \langle n+1 \rangle \,, \qquad (x_0, \ldots, x_{n+1}) \mapsto (x_0, \ldots, x_j,1,x_{j+1},\ldots, x_n)\ ,\\
\Psi(\tau_n)&: \langle n+1 \rangle \rar \langle n+1 \rangle\,,\qquad (x_0,x_1, \ldots,x_n) \mapsto (x_n,x_0, x_1, \ldots, x_{n-1})\ ,\\
\Psi(\pi^k_n)&: \langle k(n+1) \rangle \rar \langle n+1 \rangle\,,\qquad x_m \mapsto
x_{\overline{m}} \ ,\\
\end{aligned} \right.
\end{equation*}
where $\, m = 0, \,1,\,\ldots, \,k(n+1) \,$ and $\, \overline{m}$ is the remainder of $m$ modulo $n+1$.

%
 \bproof[Proof of Proposition \ref{cHopfepi}]
 The category of cocommutative (DG) Hopf algebras is equivalent to the category $\Com^{\otimes}_{\ffgr}$ of strict monoidal functors $\ffgr^{\mathrm{op}} \rar \Com_k$, with an algebra $A$ corresponding to the functor $\underline{A}\,:\,\ffgr^{\mathrm{op}} \rar \Com_k\,,\,\,\langle n \rangle \mapsto A^{\otimes n}$ (see. e.g., \cite{Pir}). Now, the epicyclic module
 $\mathrm{C}^{\Psi}_{\ast}(A)$ associated to $A$ is simply given by the composition
 $$ \mathrm{C}^{\Psi}_{\ast}(A)\,:\,\begin{diagram}[small] \Delta \Psi^{\mathrm{op}} & \rTo^{\Psi} & \ffgr^{\mathrm{op}} & \rTo^{\underline{A}}& \Com_k\end{diagram}\,,\qquad [n] \mapsto A^{\otimes n+1},$$
 where $\Psi$ is the functor \eqref{epiintog} defined above. Note that by construction  of $\Psi$, the restriction $\mathrm{C}^{\Psi}_\ast(A)$ to $\Delta C^{\mathrm{op}} \subset \Delta \Psi^{\mathrm{op}}$ coincides with the cyclic bar construction associated to $A$ as an algebra.
 \eproof
 Thus, by Proposition \ref{cHopfepi}, if $A$ is a cocommutative Hopf algebra, the cyclic module $\mathrm{C}_\ast(A)$ is equipped with extra power operations (see \eqref{powerops})  given by simplicial maps:
 \begin{equation} \la{poweropsA} p^k_\ast(A)\,:\,\mathrm{C}_\ast(A) \rar \mathrm{Sd}^k[\mathrm{C}_\ast(A)]\,, \,\, k \geq 1\ .\end{equation}
To describe these maps we need first to identify their targets $ \mathrm{Sd}^k[\mathrm{C}_\ast(A)]$. Recall that for any (DG) algebra $A$ and for any (DG) $A$-bimodule $M$, one can define a simplicial (complex of) module(s) $\mathrm{C}_\ast(A,M)= \{\mathrm{C}_n(A,M)= M \otimes A^{\otimes n}\}_{n \geq 0}$ whose homology is the Hochschild homology $\HH_\ast(A,M)$. In particular, for $M=A$, we have $\mathrm{C}_\ast(A,A)=\mathrm{C}_\ast(A)$. We will use this construction for the bimodule ${}_tA^{\otimes k}$ over the algebra $A^{\otimes k}$, where the left $A^{\otimes k}$-module structure is twisted by a cyclic permutation: i.e.
 $$ (a_1 \otimes \ldots \otimes a_k) \cdot (b_1 \otimes \ldots \otimes b_k)= a_kb_1 \otimes a_1b_2 \otimes \ldots \otimes a_{k-1}b_k \ .$$
  \blemma  \la{SD}
  $(a)$ For every $k \geq 1$, there is an isomorphism of simplicial modules
  \begin{equation} \la{SDmod}  \mathrm{Sd}^k[\mathrm{C}_\ast(A)] \stackrel{\sim}{\rar} \mathrm{C}_\ast(A^{\otimes k}, {}_tA^{\otimes k}) \end{equation}
  given $($in simplicial degree $n$$)$ by `transposition of matrices':
  \begin{equation*}
  \left. \begin{aligned}
  A^{\otimes (n+1)} \otimes \stackrel{k}{\cdots} \otimes A^{\otimes (n+1)} & \rar A^{\otimes k} \otimes \stackrel{n+1}{\cdots} \otimes A^{\otimes k}\\
  \begin{pmatrix}
  a_0 & a_{n+1} & \cdots & a_{(k-1)(n+1)}\\
  a_1 & a_{n+2} & \cdots & a_{(k-1)(n+1)+1}\\
   \cdot & \cdot & \cdots & \cdot\\
   \cdot & \cdot & \cdots & \cdot\\
   \cdot & \cdot & \cdots & \cdot\\
   a_n & a_{2n+1} & \cdots & a_{k(n+1)-1}\\
   \end{pmatrix}
    & \mapsto
  \begin{pmatrix}
  a_0 & a_{1} & \cdots & a_{n}\\
  a_{n+1} & a_{n+2} & \cdots & a_{2n+1}\\
   \cdot & \cdot & \cdots & \cdot\\
   \cdot & \cdot & \cdots & \cdot\\
   \cdot & \cdot & \cdots & \cdot\\
   a_{(k-1)(n+1)} & a_{(k-1)(n+1)+1} & \cdots & a_{k(n+1)-1}\\
   \end{pmatrix}
 \end{aligned} \right.
 \end{equation*}
 where the elements of the tensor powers $A^{\otimes (n+1)}$ and $A^{\otimes k}$ are represented as matrix columns.\\
  $(b)$ With identification \eqref{SDmod}, the power maps \eqref{poweropsA} are given by
  $$ p^k_n:A^{\otimes (n+1)} \rar A^{\otimes k} \otimes \stackrel{n+1}{\cdots} \otimes A^{\otimes k}\,,\qquad a_0 \otimes \ldots \otimes a_n \mapsto \Delta^k(a_0) \otimes \cdots \otimes \Delta^k(a_n)$$
  where $\Delta^k\,:\,A \rar A^{\otimes k}$ is the $k$-iterated the coproduct on $A$.
  \elemma
  \bproof
 $(a)$ Straightforward  verification. We leave it as an exercise to the reader.\\
 $(b)$ By definition, the maps $p^k_n$ are the images of the generating morphisms $\pi^k_n\,:\,[(n+1)k-1] \rar [n]$ of the category $\Delta\Psi$. Under the functor $\Psi$ (see \eqref{epiintog}), these morphisms correspond to the folding maps:
 $$ \nabla^k_{n+1}\,:\,\langle k(n+1) \rangle = \langle n+1 \rangle \ast \stackrel{k}{\cdots} \ast  \langle n+1 \rangle \rar \langle n+1 \rangle $$
 which act as identities $\id_{\langle n+1 \rangle}$ on each copy of the free group $\langle n+1 \rangle$ in $\langle k(n+1) \rangle$. Now, it is easy to see that the maps $\nabla^k_{n+1}$ factor in $\ffgr$ as
 $$ \begin{diagram}
   \langle n+1 \rangle \ast \stackrel{k}{\cdots} \ast \langle n+1 \rangle & \rTo^{\nabla^k_{n+1}}& \langle n+1 \rangle\\
      & \rdTo^{\cong}_{\sigma^k_n}  &  \uTo_{\nabla_k \ast \stackrel{(n+1)}{\cdots} \ast \nabla_k} \\
        & & \langle k \rangle \ast \stackrel{(n+1)}{\cdots} \ast \langle k \rangle & &\\
    \end{diagram}$$
 where $\nabla_k:\langle k \rangle = \langle 1 \rangle \ast \cdots \ast \langle 1 \rangle \rar \langle 1 \rangle$ is the $k$-folding map for $\langle 1 \rangle \in \ffgr$, and $\sigma^k_n$ is the isomorphism of free groups given by the transposition:
 $$ \begin{pmatrix}
  x_0 & x_{n+1} & \cdots & x_{(k-1)(n+1)}\\
  x_1 & x_{n+2} & \cdots & x_{(k-1)(n+1)+1}\\
   \cdot & \cdot & \cdots & \cdot\\
   \cdot & \cdot & \cdots & \cdot\\
   \cdot & \cdot & \cdots & \cdot\\
   x_n & x_{2n+1} & \cdots & x_{k(n+1)-1}\\
   \end{pmatrix}
    \mapsto
  \begin{pmatrix}
  x_0 & x_{1} & \cdots & x_{n}\\
  x_{n+1} & x_{n+2} & \cdots & x_{2n+1}\\
   \cdot & \cdot & \cdots & \cdot\\
   \cdot & \cdot & \cdots & \cdot\\
   \cdot & \cdot & \cdots & \cdot\\
   x_{(k-1)(n+1)} & x_{(k-1)(n+1)+1} & \cdots & x_{k(n+1)-1}\\
   \end{pmatrix}\,,
$$
 where the matrix columns represent the generators of the corresponding factors of the free products. Since under the functor $\underline{A}\,:\,\ffgr^{\mathrm{op}} \rar \Com_k$ the folding maps $\nabla_k:\langle k \rangle \rar \langle 1 \rangle$ correspond exactly to the $k$-iterated coproducts $\Delta^k\,:\,A \rar A^{\otimes k}$, the claim of Part $(b)$ follows.
  \eproof
 Finally, using Lemma \ref{SD}, we describe the Adams operations induced by the power maps \eqref{poweropsA} on the cyclic homology of a cocommutative DG Hopf algebra. Let $R$ be a cocommutative DG Hopf algebra which is cofibrant as an object in $\DGA_{k/k}$. Ler $R_{\n}:= R/(k+[R,R])$ denote the cyclic construction on $R$ which computes--- by a theorem of Feigin anf Tsygan (see, e.g., \cite{BKR,FT})--- the (reduced) cyclic homology $\rHC_\ast(R)$. Then, applying Lemma \ref{SD} to the epicyclic module $\mathrm{C}_\ast(R)$ we get the commutative diagram
  \begin{equation} \la{TsyganCyc} \begin{diagram}
         \mathrm{C}_{\ast}(R)   & \rTo^{p^k_\ast} & \mathrm{C}_\ast(R^{\otimes k}, {}_tR^{\otimes k})\\
             \dOnto^{\mathrm{can}}   & & \dTo_{\mu^k}\\
             R_\n & \rTo^{\overline{\Psi}^k} & R_\n\\
             \end{diagram}
              \end{equation}
 where $\mathrm{can}$ is the canonical projection onto $\pi_0 \mathrm{C}_\ast(R) \cong R_\n$ and $\mu^k$ is the composition of the natural map $\, \pi_0:\,
\mathrm{C}_\ast(R^{\otimes k}, {}_tR^{\otimes k}) \twoheadrightarrow  \big({}_tR^{\otimes k}\big)_\n \,$
 with the map $({}_tR^{\otimes k})_\n \stackrel{\overline{\mu}}{\rar} R_\n$ induced by iterated multiplication on $R$.  It follows from  Lemma \ref{SD}$(b)$ that the maps
 $\overline{\Psi}^k\,:\,R_\n \rar R_\n$ in \eqref{TsyganCyc} are induced by the compositions
 \begin{equation} \la{AdamsCoprod} \Psi^k\,:\, \begin{diagram} R & \rTo^{\Delta^k} & R^{\otimes k} & \rTo^{\mu} & R \end{diagram}\,,\,\, k \geq 1 \ .\end{equation}
  Thus, we conclude
  \bcor \la{epiCoprod}
  For any cocommutative DG Hopf algebra $R$ which is cofibrant in $\DGA_{k/k}$, the Adams operations on $\rHC_\ast(R)$ coming from the epicyclic structure on $\mathrm{C}_\ast(R)$ are induced by the maps \eqref{AdamsCoprod}.
  \ecor
  In the next section, we will give a different construction of these Adams operations in terms of derived functors, following \cite{BFPRW}.
\subsection{Lie-Hodge decomposition}
Given a Lie algebra $ \mfa $ over $k$, we consider
the symmetric ad-invariant $k$-multilinear forms on $ \mfa \,$ of a (fixed) degree $ p \ge 1 $. Every such form is induced from the universal
one: $\,\mfa \times \mfa \times \ldots \times \mfa \to \lambda^{(p)}(\mfa) \,$, which takes its values in the space $\,\lambda^{(p)}(\mfa)\,$ of coinvariants of the adjoint representation of $ \mfa $ in $ \Sym^p(\mfa)\,$.
The assignment $\,\mfa \mapsto \lambda^{(p)}(\mfa)\,$ defines a (non-additive) functor on the category of Lie algebras that naturally extends to the category of DG Lie algebras:
\begin{equation}
\la{lam}
\lambda^{(p)}:\,\DGL_k \rar \Com_k \ ,\quad \mfa \mapsto \Sym^p(\mfa)/[\mfa, \Sym^p(\mfa)]\ .
\end{equation}
The functor \eqref{lam} does {\it not} preserve quasi-isomorphisms and hence does not descend to the homotopy category
$ \Ho(\DGL_k) $. To remedy this problem, we replace $\,\lambda^{(p)}\,$ by its (left) derived functor
\begin{equation}
\la{Llam}
\L\lambda^{(p)}:\,\Ho(\DGL_k) \to \D(k)\ ,
\end{equation}
which takes its values in the derived category $ \D(k)  $
of $k$-complexes. We write  $\,\HC^{(p)}_{\bullet}(\mfa)\,$ for the homology of $\, \L\lambda^{(p)}(\mfa) \,$ and call it the {\it Lie-Hodge homology} of $ \mfa $.

For $ p = 1 $, the functor $ \lambda^{(1)} $ is just abelianization of
Lie algebras; in this case, the existence of $\, \L\lambda^{(1)} \,$ follows from
Quillen's general theory (see \cite[Chapter~II, \S 5]{Q1}), and
$\, \HC^{(1)}_{\bullet}(\mfa) \,$ coincides (up to shift in degree)
with the classical Chevalley-Eilenberg homology $ \H_\bullet(\mfa, k) $ of the Lie algebra $ \mfa $.
For $p=2$, the functor $ \lambda^{(2)} $ was introduced by Drinfeld \cite{Dr}; the existence of $ \L\lambda^{(2)} $ was established by Getzler and Kapranov \cite{GK} who
suggested that $
\HC^{(2)}_{\bullet}(\mfa) $ should be viewed as an (operadic) version of cyclic homology for Lie algebras. The existence of $\L\lambda^{(p)}$ for arbitrary $p$ was established in \cite[Sect. 7]{BFPRW}.

Next, consider the functor
$$ (\mbox{--})_{\n}\,:\,\DGA_{k/k} \rar \Com_k\,\qquad R \mapsto R/(k+[R,R])\,,$$
which is called the cyclic functor on associative DG algebras ({\it cf.} \cite{FT}). Observe that each $ \lambda^{(p)} $ comes together with a natural transformation  to the composite functor
$ \U_\n := (\,\mbox{--}\,)_\n \circ \,\U:\,\DGL_k \to \DGA_{k/k} \to \Com_k $, where
$\U $
denotes the universal enveloping algebra functor on the category of (DG) Lie algebras.  The natural transformations $\,\lambda^{(p)} \to \U_\n \,$ are induced by the symmetrization maps
\begin{equation}
\la{symfun}
\Sym^p(\mfa) \to \U\mfa\ ,\quad x_1 x_2 \ldots x_p\, \mapsto\,
\frac{1}{p!}\,\sum_{\sigma \in {\mathbb S}_p}\, \pm \,x_{\sigma(1)} \cdot x_{\sigma(2)} \cdot
\ldots \cdot x_{\sigma(p)}\ ,
\end{equation}
which, by the Poincar\'e-Birkhoff-Witt Theorem, assemble to an isomorphism of  DG $\mfa$-modules \\
$\,\Sym_k(\mfa) \cong \U \mfa \,$. From this, it follows that   $\,\lambda^{(p)} \to \U_\n \,$ give an isomorphism of functors
\begin{equation}
\la{eqv1}
\bigoplus_{p=1}^{\infty} \lambda^{(p)} \,\cong \, \U_\n \ .
\end{equation}
On the other hand, by a theorem of Feigin and Tsygan
\cite{FT} (see also \cite{BKR}), the functor $ (\,\mbox{--}\,)_\n $ has a left
derived functor $ \L(\,\mbox{--}\,)_\n:\, \Ho(\DGA_{k/k}) \to \D(k) $ that computes
the reduced cyclic homology $\,\rHC_\bullet(R)\,$ of an associative algebra
$ R \in \DGA_{k/k} $. Since $ \U $ preserves quasi-isomorphisms and maps cofibrant
DG Lie algebras to cofibrant DG associative algebras, the isomorphism \eqref{eqv1}
induces an isomorphism of derived functors from $ \Ho(\DGL_k) $ to $ \D(k) $:
\begin{equation}
\la{eqv2}
\bigoplus_{p=1}^{\infty}\, \L\lambda^{(p)}\, \cong\, \L(\,\mbox{--}\,)_\n \circ \,\U \ .
\end{equation}
At the level of homology, \eqref{eqv2} yields the direct decomposition ({\it cf.}~\cite[Theorem 7.2]{BFPRW}.

\begin{equation}
\la{hodgeds}
\rHC_{\bullet}(\U\mfa) \,\cong\,\bigoplus_{p=1}^{\infty}\, \HC^{(p)}_{\bullet}(\mfa)\ \text{.}
\end{equation}
To state the main theorem of this section, we recall that the universal enveloping algebra $\U\mfa$ of a (DG) Lie algebra $\mfa$ has the natural structure of a cocommutative (DG) Hopf algebra. By Proposition \ref{cHopfepi}, the associated simplicial module $\mathrm{C}_\ast(\U\mfa)$ carries therefore an epicyclic structure. We write $\overline{\Psi}^k\,:\,\rHC_{\ast}(\U\mfa) \rar \rHC_{\ast}(\U\mfa)$, $k \geq 1$, for the Adams operations induced by this structure.
\bthm \la{twoAdams}
For every $p \geq 2$, the Lie-Hodge homology $\HC^{(p)}_\ast(\mfa)$  is the common (graded) eigenspace of the operators $\overline{\Psi}^k$ corresponding to the eigenvalues $k^p$, $k \geq 1$.
\ethm
\bproof
Without loss of generality, we may assume that $\mfa$ is a cofibrant DG Lie algebra. Then, as explained above, $R=\U\mfa$ is a cocommutative DG Hopf algebra which is cofibrant as an associative DG algebra. By Corollary \ref{epiCoprod}, the Adams operations are induced by the maps $\mu_k \circ \Delta^k\,:\, R \rar R$, where $\Delta^k:R \rar R^{\otimes k}$ is the $k$-iterated coproduct on $R$ and $\mu_k: R^{\otimes k} \rar R$ is the $k$-iterated product. Theorem follows now from \cite[Proposition 7.5 and Corollary 7.7]{BFPRW}, which shows that the very same Adams operations $\overline{\Psi}^k$ arise from the derived functors $\L\lambda^{(p)}$.
\eproof
The above theorem implies that the Lie-Hodge decomposition \eqref{hodgeds} arises from the natural epicyclic structure on $\mathrm{C}_\ast(\U\mfa)$ given in Proposition \ref{cHopfepi}.

\subsection{Hodge decomposition of Hochschild homology} \la{secHochHodge}
The decomposition \eqref{hodgeds} also extends to (reduced) Hochschild homology (see \cite[Sect. 2.1]{BRZ}):
$$ \overline{\HH}_\bullet(\U\mfa)\,\cong\,\bigoplus_{p=0}^{\infty} \HH^{(p)}_\bullet(\mfa)\,\text{.}$$
Recall that there is a natural isomorphism $\overline{\HH}_\bullet(\U\mfa)\,\cong\,\H_\bullet(\mfa; \Sym(\mfa))$ (see \cite[Theorem 3.3.2]{L}): under this isomorphism, the summand $\HH^{(p)}(\mfa)$ is identified with $\H_\bullet(\mfa;\Sym^p(\mfa))$. The Connes periodity sequence for $\U\mfa$ decomposes into a direct sum of Hodge components: the summand of Hodge degree $p$ is given by the long exact sequence  (see \cite[Theorem 2.2]{BRZ})
\begin{equation} \la{conneshodgep}  \begin{diagram}[small] \ldots & \rTo^S & \HC_{n-1}^{(p+1)}(\mfa) & \rTo^{B} & {\HH}_n^{(p)}(\mfa) & \rTo^I & \HC_n^{(p)}(\mfa) & \rTo^S & \HC_{n-2}^{(p+1)}(\mfa) & \rTo & \ldots \end{diagram} \, \text{.} \end{equation}

Next, we shall show that the Hochschild cohomology $\HH^{\ast}(\U\mfa,\U\mfa)$ has a similar Hodge decomposition. Recall that $A\,\in\,\DGA_{k/k}$ is Koszul dual to $C\,\in\,\DGC_{k/k}$ if there is a quasi-isomorphism of DG algebras $R:=\cb(C) \stackrel{\sim}{\rar} A$, where $\cb(C)$ denotes the (associative) cobar construction of $C$. Assume that $A$ is Koszul dual to $C$. Let $\iota\,:\,C \rar R$ denote the universal twisting cochain. Further recall that given a twisting $\tau\,:\,C \rar A$, there is a convolution algebra $\Hom^{\tau}(C,A)$ with twisted differential $d_{\Hom(C,A)}+[\tau,\mbox{--}]$.
The following proposition is well known (see \cite[Theorem 1.1]{N} for instance).
\bprop \la{HHCup} There is an isomorphism of graded $k$-algebras $ \HH^{\ast}(A,A)\,\cong\, \H_{-\ast}[\Hom^{\iota}(C,R)]$.
\eprop
Recall that for a DG algebra $E$, the Hochschild cochain algebra $\mathrm{C}^{\ast}(E,E)$ (with product given by the cup product) is isomorphic to $\Pi_{n \geq 0} \Hom(\bar{E}[1]^{\otimes n}, E)$ as a graded algebra. The graded subspace $\oplus_{n \geq 0} \Hom(\bar{E}[1]^{\otimes n}, E)$ is a (DG) Gerstenhaber subalgebra of $\mathrm{C}^{\ast}(E,E)$, which we shall denote by $\mathrm{C}^{\ast}_{\oplus}(E,E)$. Let $\HH^{\ast}_\oplus(E,E)$ denote the corresponding cohomology.
\bcor \la{KosDualCochain}
If $C$ is finite dimensional, and $E:=C^{\ast}$ is the graded linear dual of $C$, then there is an isomorphism of algebras
$$ \HH^{\ast}(A,A)\,\cong\,\HH^{\ast}_{\oplus}(E,E)\ .$$
\ecor
\bproof
Since $C$ is finite dimensional,
$$\Hom^\iota(C,R)\,\cong\,\oplus_{n \geq 0} \Hom(C,\bar{C}[-1]^{\otimes n})\,\cong\, \oplus_{n \geq 0}\Hom(\bar{E}[1]^{\otimes n}, E)\ .$$
 Thus $\Hom^{\iota}(C,R)$ (viewed as a cochain complex by inverting degrees) is isomorphic to (the normalized) Hochschild cochain complex $\mathrm{C}^{\ast}_{\oplus}(E,E)$ as a graded vector space. Is is easy to verify that the above identification turns the differential on $\Hom^\iota(C,R)$ into the Hochschild differential on $\mathrm{C}^{\ast}_{\oplus}(E,E)$ and the convolution product on $\Hom^\iota(C,R)$ into the cup product on $\mathrm{C}^{\ast}_{\oplus}(E,E)$. By Proposition \ref{HHCup}, there is an isomorphism of associative algebras $\HH^{\ast}(A,A)\,\cong\,\HH^{\ast}_{\oplus}(E,E)$.
\eproof

Let $\mfa\,\in\,\DGL_k$ be Koszul dual to $C\,\in\,\cDGC_{k/k}$. This is equivalent to the existence of a quasi-isomorphism $\mathcal{L}:=\cb_{\mathtt{Comm}}(C) \stackrel{\sim}{\rar} \mfa$, where $\cb_{\mathtt{Comm}}(C)$ denotes the (Lie) cobar construction of $C$. Let $R\,:=\,\cb{C}$ and let $\iota\,:\, C \rar R$ be the universal twisting cochain. It is easy to verify that $R\,\cong\,\U\mathcal{L}$. Since the image of $\iota$ lies in $\mathcal{L} \subset R$, we may view $\iota$ as a twisting cochain from $C$ to $\mathcal{L}$ as well. Let $R^{(p)}$ denote the image of $\Sym^p(\mathcal{L})$ in $R$ under the symmetrization map \eqref{symfun}. The adjoint action of $\mathcal{L}$ on $R$ induces an action of $\Hom(C,\mathcal{L})$ on $\Hom(C,R)$: indeed, viewing $\mathcal{L}$ as a Lie subalgebra of $R$, we can consider $[\alpha, f]\,\in\,\Hom(C,R)$ for $\alpha\,\in\,\Hom(C,\mathcal{L}),f\,\in\,\Hom(C,R)$. This action equips $\Hom^{\iota}(C,R)$ with the structure of a Lie module over the DG Lie algebra $\Hom^{\iota}(C,\mathcal{L})$. Further assume that $C$ is finite dimensional.
\bthm \la{dualHodgeCoh}
The natural inclusions $R^{(p)} \hookrightarrow R$ induce a direct sum decomposition of DG \\ $\Hom^{\iota}(C,\mathcal{L})$-modules,
$$\Hom^\iota(C,R)\,\cong\,\bigoplus_{p=0}^{\infty} \Hom^{\iota}(C,R^{(p)})\ . $$
As a consequence,
$$ \HH^{\ast}(\U\mfa,\U\mfa)\,\cong\,\bigoplus_{p=0}^{\infty} \H^{\ast}(\mfa;\Sym^p(\mfa))\ .$$
\ethm
\bproof
There is an isomorphism of $\mathcal{L}$-modules $R\,\cong\,\oplus_{p=0}^{\infty} R^{(p)}$. Since $C$ is finite dimensional,
\begin{equation} \la{HodgeHHC} \Hom(C, R)\,\cong\, \bigoplus_{p=0}^{\infty} \Hom(C,R^{(p)})\end{equation}
 as graded vector spaces. It remains to check that if $\alpha\,\in\,\Hom(C,\mathcal{L})$ and if $f\,\in\,\Hom(C,R^{(p)})$, then $[\alpha,f]\,\in\,\Hom(C,R^{(p)})$. Indeed,
 \begin{eqnarray*}
 [\alpha, f](c) &=& (-1)^{|c'||f|}\alpha(c')f(c'') - (-1)^{|f||\alpha|+|\alpha||c'|}f(c')\alpha(c'')\\
                &=& (-1)^{|c'||f|}\alpha(c')f(c'')-(-1)^{|f||\alpha|+|\alpha||c''|+|c'||c''|}f(c'')\alpha(c')\\
                &=& (-1)^{|c'||f|}[\alpha(c'),f(c'')]\,,
 \end{eqnarray*}
 where the second equality above follows from the fact that $C$ is cocommutative. Since $\alpha(c')\,\in\,\mathcal{L}$ and $f(c'')\,\in\,R^{(p)}$, $[\alpha(c'),f(c'')]\,\in\,R^{(p)}$ because \eqref{symfun} is a morphism of $\mathcal{L}$ modules (with the adjoint action). This shows that \eqref{HodgeHHC} is a morphism of graded $\Hom(C,\mathcal{L})$ Lie modules. In particular, for $f\,\in\,\Hom(C,R^{(p)})$, $[\iota,f]\,\in\Hom(C,R^{(p)})$. The differential on $\Hom^\iota(C,R)$ thus restricts to $\Hom^\iota(C,R^{(p)})$ for each $p$. Hence, \eqref{HodgeHHC} is an isomorphism of complexes
 $$\Hom^\iota(C,R)\,\cong\,\bigoplus_{p=0}^{\infty}  \Hom^\iota(C,R^{(p)})\ .$$
 Next, the Jacobi identity for the commutator bracket on the convolution algebra $\Hom(C,R)$ implies that the action of $\Hom(C,\mathcal{L})$ on $\Hom(C,R)$ is compatible with the twisted differential $\partial+[\iota,\mbox{--}]$. This proves the first statement of the desired theorem. The second statement follows from Proposition \ref{HHCup} once we verify that
 \begin{equation} \la{LieCoh} \H_{-\ast}[\Hom^\iota(C,R^{(p)})]\,\cong\,\H^{\ast}(\mfa,\Sym^p(\mfa))\ . \end{equation}
 Since $C$ is Koszul dual to $R$, $R \otimes_\iota C$ is a semi-free resolution of $k$ as a DG left $R\,\cong\,\U\mathcal{L}$-module. Similarly, since $C$ is Koszul dual to $\mfa$, $\U\mfa \otimes_\tau C$ is a semi-free resolution of $k$ as a DG left $\U\mfa$-module, where $\tau$ denotes the composite twisting cochain $C \stackrel{\iota}{\rar} \mathcal{L} \stackrel{\sim}{\rar} \mfa$. It follows that there are isomorphisms in the derived category of complexes of $k$-vector spaces
 \begin{eqnarray*}
 \Hom^\iota(C,R^{(p)}) &\cong &  \Hom_{\U\mathcal{L}}(R \otimes_{\iota} C, \Sym^p(\mathcal{L}))\\
                              & \cong & \RHom_{\U\mathcal{L}}(k, \Sym^p(\mathcal{L}))\\
                                   & \cong & \RHom_{\U\mathcal{L}}(k,\Sym^p(\mfa))\\
                                   & \cong & \Hom_{\U\mathcal{L}}(R \otimes_{\iota} C, \Sym^p(\mfa))\\
                                   &\cong & \Hom^\tau(C,\Sym^p(\mfa))\\
                                   & \cong & \Hom_{\U\mfa}(\U\mfa \otimes_\tau C, \Sym^p(\mfa))\\
                                   & \cong & \RHom_{\U\mfa}(k,\Sym^p(\mfa))\ .
 \end{eqnarray*}
This implies the isomorphism \eqref{LieCoh} on homologies.
 \eproof
Let $E:=C^{\ast}$ denote the (graded) linear dual of $C$. Then, $\Hom(\bar{E}[1],E)$ is an $E$-module via the action $(x \cdot f)(y)=x \cdot f(y)$ for $f\,\in\,\Hom(\bar{E}[1],E)$, $x \,\in\,E\,,y\,\in\,\bar{E}[1]$. It is easy to verify that
$$ \Hom(C,R^{(p)})\,\cong\,\Sym^p_{E}[\mathcal{L}_{E}(\Hom(\bar{E}[1],E))]\,,$$
where $\mathcal{L}_E(V)$ denotes the free Lie algebra generated (over $E$) by a free $E$-module $V$. From Theorem \ref{dualHodgeCoh} and Corollary \ref{KosDualCochain}, it follows that
\bcor \la{HodgeHHE}
There is a direct sum decomposition
$$ \mathrm{C}^{\ast}_{\oplus}(E,E)\,\cong\,\bigoplus_{p=0}^{\infty} \Sym^p_{E}[\mathcal{L}_{E}(\Hom(\bar{E}[1],E))] \ .$$
Moreover,
$$ \H^{\ast}\big(\Sym^p_{E}[\mathcal{L}_{E}(\Hom(\bar{E}[1],E))]\big) \,\cong\,\H^{\ast}(\mfa;\Sym^p(\mfa))\ .$$
\ecor
\noindent
\textbf{Remark.} The decomposition of Hochschild cochains in Corollary \ref{HodgeHHE} is analogous to the Hodge decomposition of the complex of polydifferential operators on a smooth proper variety (over a field of characteristic $0$) in \cite{R} (see {\it loc. cit.}, Section 4).

\section{Cyclic pairings and Poisson structures} \la{seccyclicpoiss}

\subsection{Derived Poisson structures} The notion of a derived Poisson algebra was introduced in \cite{BCER} (see also \cite{BRZ}), as a natural homological generalization of noncommutative Poisson algebras in the sense of Crawley-Boevey \cite{CB}. We briefly recall basic definitions.

\subsubsection{Definitions}
\la{Defs}
Let $ A$ be an (augmented) DG algebra. The space  $\DER(A)$  of graded $k$-linear derivations of $A$
is naturally a DG Lie algebra with respect to the commutator bracket. Let $\DER(A)^\n $ denote the subcomplex of  $\DER(A) $
comprising derivations with image in $\,k+[A,A] \subseteq A \,$. It is easy to see that $ \DER(A)^\n $ is a DG Lie ideal of $ \DER(A) $,
so that $\,\DER(A)_{\natural} := \DER(A)/\DER(A)^\n $ is a DG Lie algebra.  The natural action of $ \DER(A) $ on $A$ induces a Lie algebra
action of $ \DER(A)_\n $ on the quotient space $ A_\n := A/(k+[A,A]) $. We write $\, \varrho:  \DER(A)_\n \to \END(A_\n) \,$ for the
corresponding DG Lie algebra homomorphism.

Now,  following \cite{BCER},  we define a {\it Poisson structure} on $A$ to be  a DG Lie algebra structure on $\,A_\n \,$ such that the adjoint representation
$ \mbox{\rm ad}:\, A_\n \to \END(A_\n) $ factors through $ \varrho \,$: i.~e., there is a morphism of DG Lie algebras $\,\alpha :\, A_\n \rar \DER(A)_{\natural}\,$ such that $\,\mbox{\rm ad} = \varrho \circ \alpha \,$.  It is easy to see that if $A$ is a commutative DG algebra, then a Poisson structure on $A$ is
the same thing as a (graded) Poisson bracket on $A$. On the other hand, if $A$ is an ordinary $k$-algebra (viewed as a DG algebra), then a Poisson structure on
$A$ is precisely a ${\rm H}_0$-Poisson structure in the sense of \cite{CB}.

Let $A$ and $B$ be two Poisson DG algebras, i.e. objects of $\DGA_{k/k}$ equipped with Poisson structures.
A {\it morphism} $\,f:\, A \rar B $ of Poisson algebras is then a morphism $ f: A \to  B $ in $\DGA_{k/k} $ such that $ f_{\natural}:\, A_{\natural} \rar B_{\natural} $ is a morphism of DG Lie algebras. With this notion of morphisms, the Poisson DG algebras form a category which we denote $\mathtt{DGPA}_k $.
Note that $\mathtt{DGPA}_k $ comes  with two natural functors: the forgetful functor $ U:\, \mathtt{DGPA}_k  \to \DGA_{k/k} $ and the cyclic functor
$ (\,\mbox{--}\,)_\n : \mathtt{DGPA}_k  \to \DGL_k $. We say that a morphism  $ f $ is a {\it weak equivalence} in $ \mathtt{DGPA}_k $
if $ Uf $ is a weak equivalence in $ \DGA_{k/k} $ and $ f_\n $ is a weak equivalence in $ \DGL_k $; in other words, a weak equivalence in $\mathtt{DGPA}_k $
is a quasi-isomorphism of DG algebras,   $ f: A \to B \,$,  such that the induced map  $ f_{\n}\,:\,A_\n \rar B_\n $ is a quasi-isomorphism of DG Lie algebras.

\bprop[{\cite{BRZ}}]
The category $\mathtt{DGPA}_k$ with weak equivalences specified above is a saturated homotopical category in the sense \cite{DHKS}.
\eprop

\noindent
\textbf{Remark.} It is known (see \cite{BK}) that the categories with weak equivalences form a (closed) model category-- called the category $\underline{\mathrm{Relcat}}$ of relative categories-- which is one of several (equivalent) models of an $\infty$-category (in the sense of \cite{Lu}). Thus, although we do not know whether the category $\mathtt{DGPA}_k$ carries a model structure, we can think of it as being an $\infty$-category in the sense of \cite{Lu}.

 The above proposition allows us to define a well-behaved homotopy category of Poisson algebras
$$ \Ho(\mathtt{DGPA}_k)\,:=\, \mathtt{DGPA}_k[\mathscr{W}^{-1}]\,,$$
where $\mathscr{W}$ is the class of weak equivalences.

\vspace{1ex}

\begin{definition}[\cite{BRZ}]
A {\it derived Poisson algebra} is a cofibrant associative DG algebra $A$ equipped with a
Poisson structure, which is viewed up to weak equivalence,
i.e. as an object in $ \Ho(\mathtt{DGPA}_k) $.
\end{definition}
Since the complex $A_{\n}$ computes the (reduced) cyclic homology of a cofibrant DG algebra $A$, the (reduced) cyclic homology of a derived Poisson algebra $A$ carries a natural structure of a graded Lie algebra (see \cite[Prop. 3.3]{BRZ}).

Another important result that holds for derived Poisson algebras in
$ \Ho(\mathtt{DGPA}_k) $ and that motivates our study of these objects is the following

\begin{theorem}[{\it cf.} \cite{BCER}, Theorem~2]
\la{t3s2int}
If $A$ is a derived Poisson DG algebra, then, for any $n$,
there is a unique graded Poisson bracket on the representation homology $ \HR_\bullet(A,n)^{\GL} $, such that the derived character map
$\,
\Tr_n:\, \rHC_\bullet(A) \to \HR_\bullet(A,n)^{\GL} $
is a Lie algebra homomorphism.
\end{theorem}

\subsubsection{Necklace Lie algebras}
\la{necklace}
The simplest example of a derived Poisson algebra is the tensor algebra $A=T_kV$ generated by an even dimensional $k$-vector space $V$ equipped with a symplectic form $\langle \mbox{--},\mbox{--} \rangle\,:\, V \times V \rar V$. In this case, $A$ carries a double Poisson structure in the sense of \cite{VdB}. The double bracket
$$\{\!\{ \mbox{--},\mbox{--}\}\!\}\,:\,\bar{A} \otimes \bar{A} \rar A \otimes A $$
is given by the formula
\begin{align} \la{ndbr}
\begin{aligned}
 &\{\!\{(v_1, \ldots, v_n), (w_1, \ldots, w_m)\}\!\}\,=\, \\
 &\sum_{\stackrel{i=1,\ldots,n}{j = 1,\ldots, m}}  \langle v_i, w_j \rangle  (w_1 ,\ldots, w_{j-1}, v_{i+1} ,\ldots, v_n)\otimes (v_1, \ldots, v_{i-1}, w_{j+1} ,\ldots, w_m) \,, \end{aligned} \end{align}
where $(v_1,\ldots,v_n)$ denotes the element $v_1 \otimes \ldots \otimes v_n \,\in\, T_kV$ with $v_1,\ldots,v_n\,\in\,V$. This double bracket can be extended to $A \otimes A$ by setting $\{\!\{a,1 \}\!\}=\{\!\{ 1,a\}\!\}=0$. It induces a noncommutative Poisson structure on $A$ with Lie bracket on $A_{\n}$ given by the formula
$$\{ \bar{\alpha},\bar{\beta}\}\,=\, \overline{\mu \circ \{\!\{\alpha,\beta\}\!\}}\,, $$
where $\mu:A \otimes A \rar A$ is the multiplication map and $\bar{a}$ denotes the image of $a \in A$ under the canonical projection $A \rar A_\n$. The Lie algebra $A_\n\,=\,T_kV_{\n}$ with the above bracket is called the {\it necklace Lie algebra} (see \cite{BL,G}).

\subsection{Cyclic pairings} \la{cyclicpair}

We now describe a construction of derived Poisson structures associated with cyclic coalgebras. Recall ({\it cf.}~\cite{GK}) that a graded associative $k$-algebra is called $n$-{\it cyclic} if it carries a symmetric bilinear pairing $\langle \mbox{--},\mbox{--} \rangle\,:\, A \times A \rar k$ of degree $n$ satisfying
$$ \langle ab ,c \rangle \,=\,  \langle a, bc \rangle \,,\,\,\,\,\, \forall\,\, a,b,c\,\in\,A\,\text{.}$$
Dually, a graded coalgebra $C$ is called $n$-{\it cyclic} if it carries a symmetric bilinear pairing $ \langle \mbox{--},\mbox{--} \rangle\,:\,  {C} \times {C} \rar k$ of degree $n$ satisfying
$$ \langle v', w\rangle v'' \,=\, \pm \langle v, w''\rangle w'\,,\,\,\,\,\,\, \forall\,\,v,w\,\in\,C,$$
where $v'$ and $v''$ are the two components of the coproduct $\Delta_Cv\,=\,v' \otimes v''$ written in the Sweedler notation. Note that if $A$ is a finite dimensional graded $-n$-cyclic algebra whose cyclic pairing is non-degenerate, then $C:=\Hom_k(A,k)$ is a graded $n$-cyclic coalgebra. A DG coalgebra $C$ is $n$-cyclic if it is $n$-cyclic as a graded coalgebra and
$$ \langle du, v \rangle \pm \langle u, dv \rangle \,=\, 0\,,$$
for all homogeneous $u,v \,\in\,{C}$, i.e, if $\langle \mbox{--}, \mbox{--} \rangle \,:\, {C}[n] \otimes C[n] \rar k[n]$ is a map of complexes. {\it We say that a co-augmented DG coalgebra $C\,\in\,\DGC_{k/k}$ is $n$-cyclic if $\bar{C}$ is $n$-cyclic as a non-counital DG coalgebra}.

Assume that $C\,\in\,\DGC_{k/k}$ is equipped with a cyclic pairing of degree $n$ and let $R\,:=\,\cb(C)$ denote the (associative) cobar construction of $C$. Recall that $R\,=\, T_k(\bar{C}[-1])$  as a graded $k$-algebra. For $v_1,\ldots,v_n\,\in\,\bar{C}[-1]$, let $(v_1, \ldots ,v_n)$ denote the element $v_1 \otimes \ldots \otimes v_n$ of $R$. By~\cite[Theorem~15]{BCER}, the cyclic pairing on $C$ of degree $n$ induces a double Poisson bracket of degree $n+2$ (in the sense of~\cite{VdB})
$$\{\!\{ \mbox{--}, \mbox{--}\}\!\}\,:\, \bar{R} \otimes \bar{R} \rar {R} \otimes {R}\ .$$
This double bracket is given by the formula
\begin{align} \la{dpbr}
\begin{aligned}
 &\{\!\{(v_1, \ldots, v_n), (w_1, \ldots, w_m)\}\!\}\,=\, \\
 &\sum_{\stackrel{i=1,\ldots,n}{j = 1,\ldots, m}} \pm \langle v_i, w_j \rangle  (w_1 ,\ldots, w_{j-1}, v_{i+1} ,\ldots, v_n)\otimes (v_1, \ldots, v_{i-1}, w_{j+1} ,\ldots, w_m) \,, \end{aligned} \end{align}
generalizing \eqref{ndbr}. This double bracket can be extended to $R \otimes R$ by setting $\{\!\{r,1\}\!\}\,=\, \{\!\{1,r\}\!\}=0$. Associated to~\eqref{dpbr} is the usual bracket
\begin{equation} \la{bronr} \{ \mbox{--}, \mbox{--}\} \,:=\, \mu\,\circ\, \{\!\{ \mbox{--},\mbox{--}\}\!\}\,:\, {R} \otimes {R} \rar {R}\ .\end{equation}
 Let $\n\,:\, {R} \rar R_\n$ be the canonical projection and let $\{ \mbox{--}, \mbox{--}\}\,:\, \n \circ \{ \mbox{--}, \mbox{--}\}\,:\, {R} \otimes {R} \rar R_\n$. We recall that the bimodule ${R} \otimes {R}$ (equipped with outer $R$-bimodule structure) has a double bracket (in the sense of~\cite[Defn. 3.5]{CEEY}) given by the formula
\begin{align*}
&\{\!\{ \mbox{--}, \mbox{--}\}\!\}\,\,:\,{R} \times ({R} \otimes {R}) \rar {R} \otimes ({R} \otimes {R}) \oplus ({R} \otimes {R}) \otimes {R}\,, \\
& \{\!\{r, p \otimes q\}\!\} \,:=\,  \{\!\{r,p\}\!\} \otimes q \oplus (-1)^{|p|(|r|+n)} p \otimes \{\!\{ r,q\}\!\} \,\text{.}
\end{align*}
This double bracket restricts to a double bracket on the sub-bimodule $\Omega^1R$ of $R \otimes R$ (~\cite[Corollary~5.2]{CEEY}). Let $\{\mbox{--}, \mbox{--}\}\,:\, R \otimes \Omega^1R \rar \Omega^1R$ be the map $\mu \circ \{\!\{ \mbox{--},\mbox{--}\}\!\}$, where $\mu$ is the bimodule action map and let $\{\mbox{--}, \mbox{--}\}_\n\, :\, R \otimes \Omega^1R \rar \Omega^1R_{\n}$ denote the map $\n \circ \{ \mbox{--},\mbox{--}\}$.

 As for the necklace Lie algebra, the bracket $\{ \mbox{--},\mbox{--}\}\,:\, {R} \otimes {R} \rar R_\n$ descends to a DG $(n+2)$-Poisson structure on $R$. In particular, there is a (DG) Lie bracket $\{\mbox{--},\mbox{--}\}_{\n}$  on $R_\n$ of degree $n+2$. The restriction of~\eqref{bronr} to $\bar{R}$ induces a degree $n+2$ DG Lie module structure over $R_\n$ on $\bar{R}$ and the bracket $\{\mbox{--},\mbox{--}\}_\n\,:\,R \otimes \Omega^1R \rar \Omega^1R_\n$ induces a degree $n+2$ DG Lie module structure over $R_\n$ on $\Omega^1R_\n$ (see~\cite[Proposition~3.11]{CEEY}). On homologies, we have (see~
\cite{CEEY}, Theorem~1.1 and Theorem~1.2)
\bthm \la{liestronhom}
Let $A\,\in\,\DGA_{k/k}$ be an augmented associative algebra Koszul dual to $C\,\in\,\DGC_{k/k}$. Assume that $C$ is $n$-cyclic. Then, \\
\vspace{1ex}
$(i)$ $\rHC_{\ast}(A)$ has the structure of a graded Lie algebra (with Lie bracket of degree $n+2$).\\
\vspace{1ex}
$(ii)$ $\overline{\HH}_{\ast}(A)$ has a graded Lie module structure over $\rHC_{\ast}(A)$ of degree $n+2$.\\
\vspace{1ex}
$(iii)$ The maps $S,B$ and $I$ in the Connes periodicity sequence are homomorphisms of degree $n+2$ graded Lie modules over $\rHC_{\ast}(A)$.
\ethm
The Lie bracket of degree $n+2$ on $\rHC_{\bullet}(A)$ induced by a $(n+2)$-Poisson structure on $R_{\n}$ as above is an example of  a derived $(n+2)$-Poisson structure on $A$.

\subsubsection{} \la{convention}
\noindent
\textbf{Convention.} Since we work with algebras that are Koszul dual to $n$-cyclic coalgebras, the associated Lie algebras that we work with have Lie bracket of degree $n+2$. Similarly, all Lie modules are degree $n+2$ Lie modules. We therefore, drop the prefix ``degree $n+2$" in sections that follow. To simplify terminology, we shall refer to (derived) $(n+2)$-Poisson structures simply as (derived) Poisson structures.

\subsection{Van den Bergh duality}

Assume that $A\,\in\,\DGA_{k/k}$ is Koszul dual to $C\,\in\,\DGC_{k/k}$. Let $\tau\,:\,C \rar A$ denote the twisting cochain corresponding to the quasi-isomorphism $R \stackrel{\sim}{\rar} A$, where $R:=\cb(C)$. Further, assume that $C$ is a finite-dimensional coalgebra equipped with a cyclic pairing (of degree $-n$) which is induced by a non-degenerate cyclic pairing (of degree $n$) on the graded linear dual $E=C^{\ast}$. The pairing on $C$ induces an isomorphism (complexes) $\phi\,:\,E:=\Hom_k(C,k) \rar C[-n]$ of $k$-vector spaces whose (shifted) inverse is the linear map
$$ C\,\cong\, E[n]\,, c \mapsto \langle c, \mbox{--}\rangle\ .$$
The isomorphism $\phi\,:\,E \rar C[-n]$ induces an isomorphism of DG $R$-bimodules
\begin{equation} \la{bimoddual} \Hom_k^{\iota}(C, R^e) \rar R \otimes_{\iota} C[-n] \otimes_\iota R\,,\end{equation}
where the $R$-bimodule structure on the left is induced by the ``inner" bimodule structure on $R^e$. Identifying $\Hom^{\iota}(C,R^e)\,=\,\Hom_{R^e}(R \otimes_\iota C \otimes_\iota R,R^e)$ and noting that $R \otimes_\iota C \otimes_\iota R$ is a semifree resolution of $R$, we see that the nondegenerate cyclic pairing on $C$ induces an isomorphism in $\D^{b}(R^e)$
$$ R^{\vee} \cong R[-n]\,,$$
where $R^{\vee}$ is the (derived) bimodule dual of $R$. Taking derived tensor products over $R^e$ and homology, we obtain an isomorphism
\begin{equation} \la{vdbdual} \Psi\,:\,\HH^{\ast}(A, A) \,\cong\, \HH_{n-\ast}(A)\ .\end{equation}
The isomorphism $R^{\vee} \cong R[-n]$ induces an isomorphism $\H_\ast(R^{\vee} \otimes^{\L}_{R^e} R)\,\cong\,\H_\ast(R \otimes^{\L}_{R^e} R[-n])$. The image of the identity map on $R$ (viewed as an element of $\H_0(R^{\vee} \otimes^{\L}_{R^e} R)$ under the above isomorphism) is an element $\eta\,\in\,\HH_n(A,A)$. We recall that
\blemma[{\cite[Prop. 5.5]{deVV}}] \la{capprod}

The map $\Psi$ coincides with the map $\eta \cap \mbox{--}\,:\,\HH^{\ast}(A, A) \,\cong\, \HH_{n-\ast}(A)$.
\elemma
 Let $\{\mbox{--},\mbox{--}\}\,:\,\rHC_\ast(A) \otimes \rHC_\ast(A) \rar \rHC_\ast(A)$ denote the derived Poisson bracket on $\rHC_\ast(A)$. The following result was proven for quadratic Koszul algebras in \cite{CEEY} (see the proof of Corollary 1.5 in {\it loc. cit.}). We give below a different, more direct proof in a slightly more general context.
\bprop \la{poissoncup}
The derived Poisson bracket on $\rHC_\ast(A)$ is given by
$$ \{\alpha,\beta\}\,=\,\mathrm{I}\big[\Psi[\Psi^{-1}(B(\alpha)) \cup \Psi^{-1}(B(\beta))]\big]\,,\qquad \forall \,\,\alpha\,,\beta\,\in\,\rHC_{\ast}(A)\ .$$
\eprop
\bproof
Since $R \otimes_\iota C \otimes_\iota R$ is a semifree resolution of $R$ as an $R$-bimodule, $\HH_\ast(A)\,\cong\,\HH_\ast(R)$ can be identified with the homology of the complex $(R \otimes_\iota C \otimes_\iota R) \otimes_{R^e} R$, which is isomorphic to $R \otimes C$ as graded vector spaces. The differential on $R \otimes C$ induced by that on $R \otimes_\iota C \otimes_\iota R$ is however, twisted and differs from the differential on $R \otimes_\iota C$. We let $R \otimes_\iota C_\iota$ denote $R \otimes C$ equipped with this differential. Explicitly, for $r\,\in\,R$ and $c\,\in\,C$, we have
$$ \partial_{R \otimes_\iota C_\iota}(r \otimes c)\,=\, d_Rr \otimes c+(-1)^{|r|} r\otimes d_C c + (-1)^{|c''|(|r|+|c'|)}\tau(c'')r \otimes c' -(-1)^{|r|} r \tau(c') \otimes c'' \ .$$
On the other hand, by a theorem of Feigin and Tsygan \cite{FT} (see also \cite{BKR}), $\rHC_\bullet(A)\,\cong\, \H_\ast[R_{\n}]$. It is easy to verify that the Lie bracket on $R_\n$ is given by the composite map
$$\begin{diagram} R_\n \otimes R_\n & \rTo^{\partial \otimes \partial} & (R \otimes_\iota {C}_\iota[-1]) \otimes (R \otimes_\iota {C}_\iota[-1]) & \rTo & R \otimes R \otimes {C}[-1]^{\otimes 2} & \rTo^{\mu \otimes \langle \mbox{--},\mbox{--} \rangle} &  R &\rOnto & R_\n \end{diagram}\,,$$
where the second arrow permutes factors and $\partial\,:\, R_\n \rar R \otimes_\iota {C}_\iota[-1]$ denotes the cyclic derivative. On homology, the cyclic derivative $\partial$ induces the Connes operator $B\,:\,\rHC_\ast(A) \rar \HH_{\ast+1}(A)$. It therefore, suffices to check that the map induced on homology by the composition
$$ \begin{diagram} (R \otimes_\iota {C}_\iota) \otimes (R \otimes_\iota {C}_\iota) & \rTo & R \otimes R \otimes {C}^{\otimes 2} & \rTo^{\mu \otimes \langle \mbox{--},\mbox{--} \rangle} &  R &\rOnto & R_\n \end{diagram} $$
coincides with
\begin{equation} \la{cupbracket} \begin{diagram} \HH_\ast(A) \otimes \HH_\ast(A) & \rTo^{\Psi^{-1} \otimes \Psi^{-1}} & \HH^{\ast}(A,A) \otimes \HH^{\ast}(A,A) & \rTo^{\cup}& \HH^\ast(A,A) & \rTo^{\Psi}& \HH_\ast(A,A) & \rTo^{\mathrm{I}} & \rHC_{\ast}(A)\end{diagram}\ .\end{equation}
By Proposition \ref{HHCup}, $\HH^{\ast}(A,A)$ is the homology of $\Hom^{\iota}(C,R)$, whose convolution product induces the cup product. On identifying $\Hom^\iota(C,R)$ with $R \otimes E$ as graded vector spaces, the map $\Psi^{-1}$ gets identified with the map induced on homology by $\id_R \otimes \phi^{-1}\,:\,R \otimes C \rar R \otimes E$, and the convolution product on $\Hom^{\iota}(C,R)$ is identified with the product on $R \otimes E$. On the other hand, the map $\mathrm{I}$ is induced on homology by
$$ R \otimes_{\iota} C_{\iota} \stackrel{\id_R \otimes \varepsilon_C}{\longrightarrow} R \twoheadrightarrow R_\n\ .$$
It therefore, suffices to verify that the following diagram commutes:
$$ \begin{diagram}
(R \otimes C) \otimes (R \otimes C) & \rTo & (R \otimes R) \otimes (C \otimes C) & \rTo^{\mu_R \otimes \langle \mbox{--},\mbox{--} \rangle} & R & \rTo^{\id} & R\\
    \dTo^{(\id_R \otimes \phi^{-1})^{\otimes 2}} & &   \dTo^{\id_{R^{\otimes 2}} \otimes (\phi^{-1})^{\otimes 2}} & & & & \uTo^{id_R \otimes \varepsilon_C}\\
 (R \otimes E) \otimes (R \otimes E) & \rTo & (R \otimes R) \otimes (E \otimes E) & \rTo^{\mu_R \otimes \mu_E} & R \otimes E & \rTo^{\id_R \otimes \phi}& R \otimes C\\
\end{diagram}$$
This reduces to verifying the commutativity of the diagram
$$
\begin{diagram}
C \otimes C  & \rTo^{\langle \mbox{--}, \mbox{--} \rangle} & k \\
 \uTo^{\phi \otimes \phi}  & & \uTo_{\varepsilon_C \circ \phi}\\
 E \otimes E  & \rTo^{\mu_E}& E\\
\end{diagram}
$$
Note that $\varepsilon_C\,:\,C \rar k$ coincides with $1_E$ under the identification $E=C^{\ast}$. Thus, for $v\,\in\,E$,
$$\varepsilon_C(\phi(v))\,=\,1_E(\phi(v))\,=\,\langle \phi(1_E),\phi(v) \rangle_C\,=\,\langle 1_E,v\rangle_E\,,$$
where $\langle \mbox{--},\mbox{--}\rangle_E$ denotes the original pairing on $E$. The commutativity of the above diagram therefore follows once we show that
$$ \langle \phi(v),\phi(w)\rangle\,=\,\langle v, w\rangle_E\,=\,\langle 1_E,v \cdot w\rangle_E\,,\,\,\forall\,v,w\,\in\,E\ .$$
This is a consequence of the fact that the pairing on $E$ is cyclic. This completes the proof of the desired proposition.
\eproof
Recall (see Theorem \ref{liestronhom}) that there is an action of $\rHC_\ast(A)$ on $\HH_\ast(A)$ making the latter a graded Lie module over the former. Abusing notation, we denote this action by
$$ \{\mbox{--},\mbox{--}\}\,:\,\rHC_\ast(A) \times \HH_\ast(A) \rar \HH_\ast(A)\ .$$
Let $[\mbox{--},\mbox{--}]_G$ denote the Gerstenhaber bracket on $\HH^{\ast}(A,A)$.
\bprop \la{actiononhh}
For all $\alpha\,\in\,\rHC_\ast(A)$ and for all $\beta\,\in\,\HH_\ast(A)$,
$$\{\alpha,\beta\}\,=\,\Psi\big([\Psi^{-1}(B(\alpha)),\Psi^{-1}(\beta)]_G\big)\ . $$
\eprop
\bproof
Note that the isomorphism $\Psi$ can be used to transport the Gerstenhaber bracket onto $\HH_\ast(A)$, making $\HH_\ast(A)$ a graded Lie algebra (up to shift in homological degree) with Lie bracket \\ $\Psi([\Psi^{-1}(\mbox{--}),\Psi^{-1}(\mbox{--})]_G)$. By \cite[Cor. 8.6]{CEEY} (also see {\it loc. cit.}, Proof of Theorem 1.6), $B\,:\,\rHC_\ast(A) \rar \HH_{\ast+1}(A)$ is a graded Lie algebra homomorphism, where $\rHC_\ast(A)$ is equipped with the derived Poisson bracket. It therefore, remains to verify that the action of $\rHC_\ast(A)$ on $\HH_\ast(A)$ arising out of the Lie algebra homomorphism $B$ coincides with the action arising out of the derived Poisson structure. We complete this verification in the routine computation that follows. For notational brevity, let $V:=\bar{C}[-1]$, and for $v_1,\ldots,v_n\,\in\,V$, let $(v_1,\ldots,v_n):=v_1 \otimes \ldots \otimes v_n \,\in\,R$. Pick $p=(v_1,\ldots,v_n)\,\in\,R$ and $q \otimes c=(u_1,\ldots,u_m) \otimes c\,\in\, R \otimes_{\iota} C_{\iota}$. Then,
$$ \{p, q\otimes c\}\,=\, \{p,q\} \otimes c +(-1)^{(|p|-n)|q|} \n(qd\{p,c\})\,,$$
where $\n\,:\,\Omega^1R \rar \Omega^1R_\n$ denotes the canonical projection. Hence,
\begin{eqnarray*}
\{p, q \otimes c\} &=& \sum_{\substack{1 \leqslant i \leqslant n \\ 1\leqslant k \leqslant m}} \pm \langle sv_i,su_k \rangle (u_1,\ldots,u_{k-1},v_{i+1},\ldots,v_n,v_1,\ldots,v_{i-1},u_{k+1},\ldots,u_m) \otimes c\\
&+& \n\big(\sum_{1 \leqslant j \leqslant n} \pm (u_1,\ldots,u_m)d(\langle sv_j,c\rangle v_{j+1},\ldots,v_n,v_1,\ldots,v_{j-1} ) \big)\\
&=& \sum_{\substack{1 \leqslant i \leqslant n \\ 1\leqslant k \leqslant m}} \pm \langle sv_i,su_k \rangle (u_1,\ldots,u_{k-1},v_{i+1},\ldots,v_n,v_1,\ldots,v_{i-1},u_{k+1},\ldots,u_m) \otimes c\\
&+& \sum_{1 \leqslant i < j \leqslant n} \pm \langle sv_j, c\rangle (v_{i+1},\ldots,v_{j-1},u_1,\ldots,u_m,v_{j+1},\ldots,v_n,v_1,\ldots,v_{i-1}) \otimes sv_i\\
&+&\sum_{1 \leqslant j <i \leqslant n} \pm \langle sv_j,c\rangle (v_{i+1},\ldots,v_n,v_1,\ldots,v_{j-1},u_1,\ldots,u_m,v_{j+1},\ldots,v_{i-1}) \otimes sv_i\ .
\end{eqnarray*}
On the other hand,
\begin{eqnarray*}
 \Psi^{-1}\big(B(v_1,\ldots,v_n)\big) &=& \sum_{1 \leq i \leq n} (v_{i+1},\ldots,v_n,v_1,\ldots,v_{i-1}) \otimes \widetilde{sv_i}\\
 \Psi^{-1}\big((u_1,\ldots,u_m) \otimes c\big) &=& (u_1,\ldots,u_m) \otimes \widetilde{c}\,,
 \end{eqnarray*}
 where $\widetilde{c}$ denotes $\phi^{-1}(c)$ for brevity. Hence,
\begin{eqnarray*}
[\Psi^{-1}(B(p)),\Psi^{-1}(q)]_G &=& \sum_{1 \leqslant i < j \leqslant n} \pm sv_j(\widetilde{c})(v_{i+1},\ldots,v_{j-1},u_1,\ldots,u_m,v_{j+1},\ldots,v_n,v_1,\ldots,v_{i-1}) \otimes \widetilde{sv_i}\\
&+&\sum_{1 \leqslant j <i \leqslant n} \pm sv_j(\widetilde{c}) (v_{i+1},\ldots,v_n,v_1,\ldots,v_{j-1},u_1,\ldots,u_m,v_{j+1},\ldots,v_{i-1}) \otimes \widetilde{sv_i}\\
&+& \sum_{\substack{1 \leqslant i \leqslant n \\ 1\leqslant k \leqslant m}} \pm sv_i(\widetilde{su_k}) (u_1,\ldots,u_{k-1},v_{i+1},\ldots,v_n,v_1,\ldots,v_{i-1},u_{k+1},\ldots,u_m) \otimes \widetilde{c}\ .
\end{eqnarray*}
Since $u(\widetilde{w})=\langle u, w\rangle$ for all $u,w\,\in\,C$, the above computations show that
$$\Psi^{-1}\big(\{p,q \otimes c\}\big)\,=\,[\Psi^{-1}(B(p)),\Psi^{-1}(q \otimes c)]_G\ . $$
Since $\rHC_\ast(A)\,=\,\H_\ast(R_\n)$ and $\HH_\ast(A)=\H_\ast(R \otimes_\iota C_\iota)$, the desired verification is complete once we apply $\Psi$ to both sides of the above equation.
\eproof
Further assume that $C$ is cocommutative, so that $A\,\cong\,\U\mfa$, where $\mfa\,\in\,\DGL_k$ is Koszul dual to $C$. The image of the counit $\varepsilon_C\,\in\,C^{\vee}$ under the isomorphism $C^{\vee} \cong C[-n]$ defines an $n$-cycle in $C$, whose class in $\H_n(C)\,\cong\,\H_n(\mfa;k)$ is denoted by $\eta$. This, in turn, defines a cap product $\eta \cap \mbox{--}\,:\,\H^{n-r}(\mfa;N) \rar \H_r(\mfa;N)$ for any DG $\mfa$-module $N$ (see \cite[Sect. 7.1]{FTV}).
\blemma \la{hodgeduality}
Under the natural isomorphisms $\HH_\ast(\U\mfa)\,\cong\,\H_\ast(\mfa;\U\mfa)$ and $\HH^\ast(\U\mfa,\U\mfa)\,\cong\, \H^\ast(\mfa;\U\mfa)$  the map $\Psi$ is identified with
$$\eta \cap \mbox{--}\,:\,\H^{\ast}(\mfa;\U\mfa) \rar \H_{n-\ast}(\mfa;U\mfa)\ . $$
As a consequence,
$$\Psi[\H^\ast(\mfa;\Sym^p(\mfa))] \,=\,\H_{n-\ast}(\mfa;\Sym^p(\mfa)) \,\forall\,p \geq 0\ .$$
\elemma
\bproof
The map $\Psi$ is induced by the map of complexes
\begin{equation} \la{vdbdual1} \Hom^{\iota}(C,R) \,\cong\, R \otimes_{\iota} C_\iota[-n] \end{equation}
obtained by tensoring the bimodule map \eqref{bimoddual} with $R$ over $R^e$. Hence, identifying the left hand side (as graded vector spaces) with $R \otimes E$, we see that \eqref{vdbdual} coincides with the map $\id_R \otimes \phi$.

On the other hand, the natural isomorphisms  $\HH_\ast(\U\mfa)\,\cong\,\H_\ast(\mfa;\U\mfa)$ and $\HH^\ast(\U\mfa,\U\mfa)\,\cong\, \H^\ast(\mfa;\U\mfa)$  are induced by the maps $\Hom^{\iota}(C, R) \rar \Hom^{\tau}(C,\U\mfa)$ and $R \otimes_\iota C_\iota \rar \U\mfa \otimes_\tau C_\tau$ induced by the canonical projection $R \stackrel{\sim}{\rar} \U\mfa$, and where $\tau:C \rar \U\mfa$ denotes the twisting cochain corresponding to the algebra homomorphism $R \stackrel{\sim}{\rar} \U\mfa$. We therefore, need to verify that the map $\eta \cap \mbox{--}$ is induced on homologies by the map
$$ \id_{\U\mfa} \otimes \phi\,:\,\Hom^{\tau}(C,\U\mfa) \rar \U\mfa \otimes_\tau C_\tau[-n]\ .$$
Identifying $\Hom^{\tau}(C,\U\mfa)\,=\,\U\mfa \otimes E$ as a graded vector space, we see that by \cite[Sect. 7.1]{FTV}, the map $\eta \cap \mbox{--}$ is induced on homologies by a map of complexes which coincides (as a map of graded vector spaces) with
$$\begin{diagram}[small] \U\mfa \otimes E  & \rTo^{\id_{\U\mfa \otimes E} \otimes \eta} &  \U\mfa \otimes E \otimes C[-n] & \rTo^{\id \otimes\Delta}& \U\mfa \otimes  E \otimes C \otimes C[-n] & \rTo^{\id_{\U\mfa} \otimes \mathrm{ev} \otimes \id} \U\mfa \otimes C[-n] \ . \end{diagram}$$
It therefore, suffices to check that $\Phi\,:\,E \rar C[-n]$ coincides with the map
$$ \begin{diagram} E  & \rTo^{\id_E \otimes \eta} & E \otimes C[-n] & \rTo^{id_E \otimes\Delta}& E \otimes C \otimes C[-n] & \rTo^{\mathrm{ev} \otimes \id} C[-n]\end{diagram}$$
which is clear. This completes the proof of Lemma \ref{hodgeduality}.
\eproof
\section{Hodge decomposition of derived Poisson structures}

\la{SecDerivedPoissonHodge}

\subsection{The main theorem}
\la{minsullivan}
Recall that if $ \g $ is an $L_\infty$-algebra with higher operations (Lie brackets)
$\, m_k: \wedge^k \g \to \g$, $\, k \ge 1 $, the lower central filtration of $ \g $ is
defined inductively by
\begin{equation}
\la{lowseries}
F^1\g := \g\ ,\quad  F^r\g \, :=\, \sum_{i_1+\ldots+i_k=r} m_k(F^{i_1}\g,\ldots,F^{i_k}\g)\, , \ r \ge 2\ .
\end{equation}
Then $ \g $ is called {\it nilpotent} if its lower central filtration \eqref{lowseries} terminates
after finitely many steps, i.e. $F^r\g=0$ for $r \gg 0$ (see, e.g., \cite[Definition 4.2]{Ge}).
For the rest of this section, we make the following
\vspace*{1ex}

\noindent
{\bf Assumption.} {\it  $\g$ is a non-negatively graded, finite-dimensional, nilpotent $L_\infty$-algebra}.

\vspace*{1ex}

Let $\mathcal{A}=(\Sym(V),Q)$ be the Chevalley-Eilenberg cochain algebra of $\g$, where $V=\g^{\ast}[-1]$, and let $W:=\g[1]$ denote the (graded) linear dual of $V$. Following \cite[Sect.5]{CR}, we will use the language of formal differential geometry, regarding $W$ as a supermanifold and $\mathcal{A}$ as the algebra of functions on $W$ equipped with cohomological vector field $Q$ of (cohomological) degree $1$. Note that the algebra $\Sym(V) \otimes \Sym(W[-1])$ of polyvector fields on $W$ is naturally a graded Lie module over the Lie algebra $\Der(\Sym(V))\,\cong\,\Sym(V) \otimes W$ of (graded) derivations of $\Sym(V)$. Let $\mathcal{V}$ denote the algebra of polyvector fields on $W$ equipped with the differential given by the action of the derivation $Q\,\in\,\Der(\Sym(V))$. The Schouten bracket makes $\mathcal{V}$ a Gerstenhaber algebra. The Hochschild-Kostant-Rosenberg map $\mathrm{I}_{\rm HKR}$ maps $\mathcal{V}$ to the complex $\mathrm{D}^{\ast}_{\mathrm{poly}}(\mathcal{A})$, which is a subcomplex of $\mathrm{C}^{\ast}_{\oplus}(\mathcal{A},\mathcal{A})$ spanned by multilinear maps $\mathcal{A}^{\otimes n} \rar \mathcal{A}$ that are differential operators in each argument. Abusing notation, we continue to denote the composite map
$$\begin{diagram} \mathcal{V} &\rTo^{\mathrm{I}_{\rm HKR}} & \mathrm{D}^{\ast}_{\mathrm{poly}}(\mathcal{A}) & \rInto & \mathrm{C}^{\ast}_{\oplus}(\mathcal{A},\mathcal{A}) \end{diagram} $$
by $\mathrm{I}_{\rm HKR}$.
The following proposition is a consequence of the Duflo-Kontsevich Isomorphism Theorem for symmetric algebras equipped with (co)homological differential (see, e.g., \cite[Theorem 5.3]{CR}):
\bprop \la{DufloRE}
The map $\mathrm{I}_{\rm HKR}\,:\,\mathcal{V} \rar  \mathrm{C}^{\ast}_{\oplus}(\mathcal{A},\mathcal{A})$ is a quasi-isomorphism of complexes that induces an isomorphism of algebras on cohomology.
\eprop
\bproof
Recall that the graded algebra of differential forms on $W$ is defined by $\Omega(W)\,:=\,\Sym(V \oplus V[1])$. For any $x\,\in\,V$ we write $dx$ for the corresponding element in $V[1]$. The de Rham differential is the derivation of (homological) degree $1$ given on generators by $d(x)=dx,\, d(dx)=0$. There is an action $\iota$ of differential forms on polyvector fields by contraction, where $x\,\in\,V$ acts by left multiplication and $dx$ acts by the derivation $\iota_{dx}$ such that $\iota_{dx}(y)=0$ for $y\,\in\,V$ and $\iota_{dx}(s^{-1}v)=x(v)$ for $v\,\in\,W$. Choosing a basis $\{x^1,\ldots,x^n\}$ in $V$ that
consists of homogeneous elements, we define the $1$-form $\alpha\,\in\,\Omega^1(W) \otimes \End(W)$ whose matrix with respect to the basis $\{e_i:=\partial_{x^i}\}$ of $W$ is given by
$$ \alpha^j_i\,:=\, d(\partial_{x^i}Q(x^j))\,=\, \frac{\partial^2Q(x^j)}{\partial x^i \partial x^k} dx^k\ .   $$
Let  $\widehat{\Omega}(W)$ denote the completion of $\Omega(W)$ with respect to the ideal generated by $V[1]$. The {\it Todd genus} associated to $ \alpha $ is defined by
$$J(\alpha):= \mathrm{Ber}\left(\frac{1- \mathrm{exp}(-\alpha)}{\alpha}\right)\,\in\,\widehat{\Omega}(W)\, ,$$
where $\mathrm{Ber}:\widehat{\Omega}(W) \otimes \End(W) \rar \widehat{\Omega}(W)$ is a map induced by the Berezinian on $\End(W)$. Note that $\frac{1-\mathrm{exp}(-\alpha)}{\alpha}\,=\,\mathrm{exp}(\sum_k c_k\alpha^k)$ for some formal power series $\sum_k c_k \alpha^k$ with constant term $1$. Thus,
$$J(\alpha)\,=\,\mathrm{exp}\big(\sum_{k=1}^{\infty} c_k\mathrm{Str}(\alpha^k)\big)\,,$$
where $\mathrm{Str}\,:\,\Omega(W) \otimes \End(W) \rar \Omega(W)$ is a linear map induced by the (super)trace on $\End(W)$.

We claim that $J(\alpha)=1$. It suffices to show that $\alpha$ is nilpotent, since then $\mathrm{Str}(\alpha^k)=0$ for all $k$. The nilpotency of $\alpha$ follows from the nilpotency of the  $L_{\infty}$-algebra $\g=W[-1]$. Recall that the (co)restriction of the Chevalley-Eilenberg (co)derivation on $\C^{\ast}(\g;k)$ to $W=\g[1]$ gives a linear map $\Sym(W) \rar W$, whose component in degree $p$ we denote by
$$ [\mbox{--},\ldots,\mbox{--}]_p\,:\,\Sym^p(W) \rar W\,,\qquad v_1 \cdot \ldots \cdot v_p \mapsto [v_1,\ldots,v_p]_p\ .$$
Note that $[v_1,\ldots,v_p]_p\,=\,\pm m_p(s^{-1}v_1,\ldots,s^{-1}v_p)$ for $v_1,\ldots,v_p$ homogeneous in $W$, where $m_p: \g^{\otimes p} \rar \g$ is the $p$-fold bracket on $\g$ (which, of course, is antisymmetric).
Now,
\begin{equation} \la{qxj} Q(x^j)= -\sum_{p \geqslant 1; i_1,\ldots,i_p} \frac{1}{p!}[e_{i_1},\ldots,e_{i_p}]^j_p x^{i_1} \ldots x^{i_p}\,,\end{equation}
where $v^j$ stands for the coefficient of $e_j$ in $V$ for any $v \in W$. Therefore,
\begin{eqnarray*}
\alpha(e_i)\,=\, \partial_{x_i}dQ(x^j)e_j &= & -\sum_{p \geqslant 2; i_1,\ldots,i_p} \frac{1}{(p-2)!}[e_{i_1},e_{i_2},\ldots,e_{i_p}]_p dx^{i_1}\delta_{ii_2} \ldots x^{i_p}\\
              &=& -\sum_{p \geqslant 2; i_1,i_3,\ldots,i_p} \frac{1}{(p-2)!}[e_{i_1},e_{i},\ldots,e_{i_p}]_p dx^{i_1}x^{i_3}\ldots x^{i_p} \,\in\, W \otimes \Omega^1(W)\,,
\end{eqnarray*}
where $\delta_{ii_2}$ is the Kronecker delta. It follows that for any (homogeneous) element $v \in W$,
\begin{equation} \la{endalpha} \alpha(v)\,=\,  \pm \sum_{p \geqslant 2; i_1,i_3,\ldots,i_p} \frac{1}{(p-2)!}[e_{i_1},v,\ldots,e_{i_p}]_p dx^{i_1}x^{i_3}\ldots x^{i_p} \,\in\, W \otimes \Omega(W)\ . \end{equation}
If $\{F^r\g\}_{r \geq 1}$ is the lower central filtration of $\g$ (see \eqref{lowseries}), we set $F^rW:=F^r\g[1]$ for $r \geq 1$. By \eqref{endalpha}, $\alpha(F^jW \otimes \Omega(W))\,\subset F^{j+1}W \otimes \Omega(W)$ for all $j$. Hence,
the image of $\alpha^r$ is contained in $F^{r+1}W \otimes \Omega(W)$ for all $r$. Since $F^rW=0$ for $r \gg 0$, $\alpha^r=0$ for $r \gg 0$. This verifies that $\alpha$ is nilpotent, as desired.

By \cite[Thm. 4.3]{LSX} (also see \cite[Thm. 5.3]{CR}), the map $\mathrm{I}_{\rm HKR} \circ \iota_{J(\alpha)^{\frac{1}{2}}}$ defines a quasi-isomorphism of complexes that induces an algebra isomorphism on cohomology. Since $J(\alpha)=1$, the desired proposition follows from Lemma \ref{Inc} below.
\eproof
\blemma \la{Inc}
The inclusion $\mathrm{D}^{\ast}_{\mathrm{poly}}(\mathcal{A}) \hookrightarrow \mathrm{C}^{\ast}_{\oplus}(\mathcal{A},\mathcal{A})$ is a quasi-isomorphism of complexes that induces an algebra isomorphism on cohomology.
\elemma
\bproof
The inclusion $\mathrm{D}^{\ast}_{\mathrm{poly}}(\mathcal{A}) \hookrightarrow \mathrm{C}^{\ast}_{\oplus}(\mathcal{A},\mathcal{A})$ is clearly compatible with cup products. It therefore, suffices to check that it is a quasi-isomorphism. Let $\mathcal{A}_0:=(\Sym(V),0)$ denote the DG algebra with trivial differential that is isomorphic to $\mathcal{A}$ as a graded algebra. First, we check that the HKR map
\begin{equation} \la{Symalg} \mathrm{I}_{\rm HKR}\,:\,\mathcal{A}_0 \otimes \Sym(W[-1]) \rar \mathrm{C}^{\ast}_{\oplus}(\mathcal{A}_0,\mathcal{A}_0) \end{equation}
is a quasi-isomorphism. Indeed, $\mathrm{C}^{\ast}_{\oplus}(\mathcal{A}_0,\mathcal{A}_0)$ may be identified with the complex $\Hom(\mathcal{A}_0 \otimes_{\pi} \bB \mathcal{A}_0 \otimes_{\pi} \mathcal{A}_0, \mathcal{A}_0)$, where $\Hom$ is in the category of {\it graded} $\mathcal{A}_0$-bimodules. Here, $\mathcal{A}_0 \otimes_{\pi} \bB \mathcal{A}_0 \otimes_{\pi} \mathcal{A}_0$ is viewed as the free resolution of $\mathcal{A}_0$ (as a graded $\mathcal{A}_0$-bimodule) whose term in homological degree $i$ is $\mathcal{A}_0 \otimes \overline{\mathcal{A}_0}^{\otimes i} \otimes \mathcal{A}_0$. On the other hand, since $\Sym^c(V[1])$ is Koszul dual to $\mathcal{A}_0$, $\mathcal{A}_0 \otimes_{\tau} \Sym^c(V[1]) \otimes_{\tau} \mathcal{A}_0$ also yields a free resolution of $\mathcal{A}_0$ as a graded $\mathcal{A}_0$-bimodule, where $\tau\,:\,\Sym^c(V[1]) \rar \mathcal{A}_0$ denotes the corresponding twisting cochain. The map $\Sym^c(V[1]) \rar \bB \mathcal{A}_0$ of coalgebras corresponding to $\tau$ induces a quasi-isomorphism $ \mathrm{C}^{\ast}_{\oplus}(\mathcal{A}_0,\mathcal{A}_0) \rar \mathcal{A}_0 \otimes \Sym(W[-1])$ whose inverse is the HKR map \eqref{Symalg}.

A standard spectral sequence argument shows that $\mathrm{I}_{\rm HKR}\,:\,\mathcal{V} \rar  \mathrm{C}^{\ast}_{\oplus}(\mathcal{A},\mathcal{A})$ is a quasi-isomorphism. Since $\mathrm{I}_{\rm HKR}$ factors through  $\mathrm{I}_{\rm HKR}\,:\,\mathcal{V} \rar \mathrm{D}^{\ast}_{\mathrm{poly}}(\mathcal{A})$, the inclusion of $\mathrm{D}^{\ast}_{\mathrm{poly}}(\mathcal{A})$ in $ \mathrm{C}^{\ast}_{\oplus}(\mathcal{A},\mathcal{A})$ is a quasi-isomorphism.
\eproof
Note that the arguments of \cite[Sect. 4]{R} go through with obvious modifications in the DG setting showing that $\mathrm{D}^{\ast}_{\mathrm{poly}}(\mathcal{A})$ splits into a direct sum of subcomplexes
$$ \mathrm{D}^{\ast}_{\mathrm{poly}}(\mathcal{A})\,\cong\, \bigoplus_{p=0}^{\infty} \Sym^p(\mathcal{L}_{\mathcal{A}}(\mathrm{D}(\mathcal{A})))\,, $$
where $\mathcal{L}_{\mathcal{A}}(\mathrm{D}(\mathcal{A}))$ denotes the free Lie algebra generated over $\mathcal{A}$ by the algebra $D(\mathcal{A})$ of differential operators on $\mathcal{A}$ (to which the Hochschild differential indeed restricts). The $\mathcal{A}$-module structure on $D(\mathcal{A})$ is given by the natural left multiplication, which makes $D(\mathcal{A})$ a semi-free DG $\mathcal{A}$-module. Setting $\mathcal{V}^p\,:=\,\mathcal{A} \otimes \Sym^p(W[-1])$, we see as in \cite[Sect. 4.2]{R} that the map $\mathrm{I}_{\rm HKR}$ restricts to a quasi-isomorphism
$$ \mathrm{I}_{\rm HKR}\,:\,\mathcal{V}^p \rar \Sym^p(\mathcal{L}_{\mathcal{A}}(\mathrm{D}(\mathcal{A}))) $$
for each $p$. Now, the subcomplex $\Sym^p(\mathcal{L}_{\mathcal{A}}(\mathrm{D}(\mathcal{A})))$ of $\mathrm{D}^{\ast}_{\mathrm{poly}}(\mathcal{A})$ is the image of a projection operator
$$e^{(p),\ast}\,:\,\mathrm{D}^{\ast}_{\mathrm{poly}}(\mathcal{A}) \rar \mathrm{D}^{\ast}_{\mathrm{poly}}(\mathcal{A})\,, $$
that is the restriction of a projection operator $e^{(p),\ast}$ on $\mathrm{C}^{\ast}_{\oplus}(\mathcal{A},\mathcal{A})$. The latter has an explicit combinatorial definition in terms of Eulerian idempotents (see \cite[Sec. 4.5]{L}). The operators $e^{(p),\ast}$ define the Hodge decomposition of Hochschild cohomology:
$$\mathrm{C}^{\ast}_{\oplus}(\mathcal{A},\mathcal{A})\,\cong\,\bigoplus_{p=0}^{\infty} \mathrm{C}^{(p),\ast}_{\oplus}(\mathcal{A},\mathcal{A})\,,\qquad \HH^\ast_{\oplus}(\mathcal{A},\mathcal{A})\,\cong\, \bigoplus_{p=0}^{\infty} \HH^{(p),\ast}_{\oplus}(\mathcal{A},\mathcal{A})\ .$$
\blemma \la{HodgeCup}
The cup product on $\HH^{\ast}_{\oplus}(\mathcal{A},\mathcal{A})$ preserves the Hodge decomposition. More precisely,
\begin{flalign*}
\HH^{(p),\ast}_{\oplus}(\mathcal{A},\mathcal{A}) \cup \HH^{(q),\ast}_{\oplus}(\mathcal{A},\mathcal{A}) \, \subseteq\,\HH^{(p+q),\ast}_{\oplus}(\mathcal{A},\mathcal{A})
\end{flalign*}
\elemma
\bproof
The previous argument shows that the Hochschild-Kostant-Rosenberg map intertwines the natural Hodge decomposition $\mathcal{V}\,\cong\,\oplus_p \mathcal{V}^p$ with that of $\mathrm{C}^{\ast}_{\oplus}(\mathcal{A},\mathcal{A})$. By Lemma \ref{Inc}, it induces isomorphisms for each $p \geq 0$
$$\mathrm{I}_{\rm HKR}\,:\, \H^\ast[\mathcal{V}^p]\,\cong\, \HH^{(p),\ast}_{\oplus}(\mathcal{A},\mathcal{A}) \ .$$
Since the product on $\mathcal{V}$ preserves the Hodge decomposition $\mathcal{V}\,\cong\,\oplus_p \mathcal{V}^p$, the desired proposition follows from Proposition \ref{DufloRE}.
\eproof


\bprop \la{cupMinLinfty}
Let $\mfa\,\in\,\DGL_k$ be as in Theorem \ref{genmain}, and let $\g$ denote the minimal $L_\infty$-model of $\mfa$. Then, the cup product and the Gerstenhaber bracket on $\HH^\ast(\U\mfa,\U\mfa)$ preserve the Lie-Hodge decomposition, i.e.,
\begin{equation} \la{cupHodgeQuillen} \HH^{(p),\ast}(\mfa) \cup \HH^{(q),\ast}(\mfa) \,\subseteq\, \,\HH^{(p+q),\ast}(\mfa)\,, \qquad
[\HH^{(p),\ast}(\mfa),\HH^{(q),\ast}(\mfa)]_G \subseteq  \HH^{(p+q-1),\ast}(\mfa)\ .\end{equation}
\eprop
\bproof
We begin by showing that there is a (Hodge decomposition preserving) isomorphism of algebras $\HH^{\ast}(\U\mfa,\U\mfa)\,\cong\,\HH^{\ast}_{\oplus}(\mathcal{A},\mathcal{A})$, where $\mathcal{A}$ is the Chevalley-Eilenberg cochain algebra of $\g$. Since $\mathcal{A}$ is neither finite-dimensional (it is only {\it locally} finite dimensional) nor bigraded, we are not in a position to quote results from the literature (see \cite[Thm. 3.5]{Ke2} and \cite{Her} for instance) for our purpose.  Let $\mathcal{C}:=\C_\ast(\g;k)$ denote the Chevalley-Eilenberg {\it chain} coalgebra of $\g$. Note that $\mathcal{C}=\mathcal{A}^{\ast}$, the graded linear dual of $\mathcal{A}$. The obvious degree $-1$ map $\mathcal{C}:=\C_\ast(\g;k) \rar \g$ is a (generalized) twisting cochain, which we denote by $\eta$. The unit of the adjunction $\cb_{\mathtt{Com}}\,:\, \cDGC_{k/k} \rightleftarrows \DGL_k\,:\,\C_\ast$ gives a weak-equivalence of DG coalgebras $\C \rar \C_\ast(\cb_{\mathtt{Com}}(\C);k)$ (see, e.g., \cite[Sec. 6.2]{BFPRW}). Let $\iota\,:\,\C \rar \cb_{\mathtt{Com}}(\C)$ denote the corresponding (canonical) twisting cochain and let $\boldsymbol{f}=(f_1,f_2,\ldots)\,:\,\g \rar \cb_{\mathtt{Com}}(\C)$ denote the corresponding $L_{\infty}$-morphism. Composition with $\boldsymbol{f}$ defines an $L_\infty$-morphism $\Hom(\C,\g) \rar \Hom(\C,\cb_{\mathtt{Com}}(\C))$, which we continue to denote by $\boldsymbol{f}$. Explicitly, for $\varphi_1,\ldots,\varphi_n\,\in\,\Hom(\C,\g)$,
$$ f_n\big(\varphi_1 ,\ldots, \varphi_n\big)\,=\, f_n \circ (\varphi_1 \otimes \ldots \otimes \varphi_n) \circ \Delta^{n-1}\,,$$
where $\Delta^{n-1}\,:\,\C \rar \C^{\otimes n}$ denotes the $n$-iterated coproduct for $n \geq 2$ (wit $\Delta^0=\id$). Clearly, $\boldsymbol{f}(\eta)=\iota$. By \cite[Prop. 1]{VD}, there is a twisted $L_\infty$-morphism
$\boldsymbol{f}^{\eta}\,:\,\Hom^{\eta}(\C,\g) \rar \Hom^{\iota}(\C,\cb_{\mathtt{Com}}(\C))$, whose $n$-th Taylor coefficient is given by the formula (see {\it loc. cit.}, Equation 2.49)
$$ f^{\eta}_n(\varphi_1 ,\ldots,\varphi_n)\,=\,\sum_{k=0}^{\infty} \frac{1}{k!} f_{n+k}(\eta ,\ldots,\eta,\varphi_1 ,\ldots,\varphi_n)\ .$$
With this, it is not difficult to verify that the following diagram commutes:
\begin{equation} \la{diffops} \begin{diagram}
\Hom^{\eta}(\C,\g) & \rTo^{f^{\eta}_1} & \Hom^{\iota}(\C,\cb_{\mathtt{Com}}(\C))\\
    \uTo^{\cong}  & & \uInto\\
  \mathcal{V}^1 & \rTo^{(\mathrm{I}_{\mathrm{HKR}})|_{\mathcal{V}^1}} & \Hom(\bar{\mathcal{A}}[1], \mathcal{A})  
\end{diagram}\ .\end{equation}
Here, the vertical arrow on the right is given by taking graded linear duals followed by composition by the inclusion $\bar{\C}[-1] \hookrightarrow \cb_{\mathtt{Com}}(\C)$. 

Since $\Sym^p(\g)$ is finite-dimensional for each $p$, there is an isomorphism of DG algebras
$$\bigoplus_p \Hom^{\eta}(\mathcal{C},\Sym^p(\g))\,\cong\, \bigoplus_p \mathcal{V}^p\,,$$
where the (convolution) algebra structure on the left hand side is induced by the coproduct on $\mathcal{C}$ and the product on $\Sym(\g)$. As in Theorem \ref{genmain}, let $C$ denote a finite dimensional cocommutative DG coalgebra Koszul dual to $\mfa$. By \cite[Thm. 10.3.15]{LV}, there is an $L_{\infty}$-quasi-isomorphism $\mfa \stackrel{\sim}{\rar} \g$. It follows that there are weak-equivalences of DG coalgebras
\begin{equation} \la{weakequcoalg} C \rar \C_\ast(\mfa;k) \rar \mathcal{C}\ .\end{equation}
Let $\tau$ denote the (generalized) twisting cochain $C \rar \mathcal{C} \stackrel{\eta}{\rar} \g$. Let $\rho$ denote the (generalized) twisting cochain $C \stackrel{\eqref{weakequcoalg}}{\rar} \C \stackrel{\iota}{\rar} \cb(\C)$. Commutativity of the diagram \eqref{diffops} implies that there is a commutative diagram
$$ \begin{diagram}
\bigoplus_p \Hom^{\eta}(\mathcal{C},\Sym^p(\g)) & \rTo^{\mathrm{I}_{\mathrm{HKR}}} & \mathrm{C}^{\ast}_{\oplus}(\mathcal{A},\mathcal{A}) & \rTo & \mathrm{C}^{\ast}_{\oplus}(\mathcal{A},E)\\
\dTo & & & & \dTo^{\cong}\\
\bigoplus_p \Hom^{\tau}(C,\Sym^p(\g)) & \rTo^{\oplus_p \Sym^p(f^{\tau}_1)} & \bigoplus_p \Hom^{\rho}(C,\Sym^p(\cb_{\mathtt{Com}}(\C))) & \rTo & \Hom^{\rho}(C,\cb(\mathcal{C}))\\
\end{diagram}\,,
$$
where the second arrow in the upper row is induced by the algebra homomorphism $\mathcal{A} \rar E:=C^{\ast}$ obtained by applying graded linear duals to \eqref{weakequcoalg}, and the last arrow on the lower row is induced by the natural symmetrization map. Note that $\Hom^{\eta}(\C,\Sym^p(\g))$ (resp., $\Hom^{\tau}(C,\Sym^p(\g))$) may be identified with $\Hom_{\DGCoMod_{\C}}(\C,\C \otimes_{\eta} \Sym^p(\g))$ (resp., $\Hom_{\DGCoMod_{\C}}(C,\C \otimes_{\eta} \Sym^p(\g))$). Under this identification, the map $\Hom^{\eta}(\C,\Sym^p(\g)) \rar \Hom^{\tau}(C,\Sym^p(\g))$ is induced by the weak-equivalence $C \rar \C$ of DG $\C$-comodules. Since $\C \otimes_{\eta} \Sym^p(\g)$ is a fibrant DG $\C$-comodule, this map is a quasi-isomorphism for each $p$. Since $\boldsymbol{f}$ is a quasi-isomorphism and since $\tau$ induces a weak-equivalence, the map $f^{\tau}_1$ is a quasi-isomorphism by \cite[Prop. 1]{VD}. Since $\mathrm{I}_{\mathrm{HKR}}$ is a quasi-isomorphism by Lemma \ref{Inc} and since $\id_\ast$ is an isomorphism of complexes, the map
$ \mathrm{C}^{\ast}_{\oplus}(\mathcal{A},\mathcal{A}) \rar \mathrm{C}^{\ast}_{\oplus}(\mathcal{A},E)$ is a quasi-isomorphism of DG algebras. Hence, there is a zig-zag of DG algebra maps
$$ \mathrm{C}^{\ast}_{\oplus}(\mathcal{A},\mathcal{A}) \rar \mathrm{C}^{\ast}_{\oplus}(\mathcal{A},E) \leftarrow  \mathrm{C}^{\ast}_{\oplus}(E,E)\,,$$
where the first arrow has been shown to be a quasi-isomorphism above and the second arrow is a quasi-isomorphism by a standard spectral sequence argument. Since the Hodge decomposition on all the above Hochschild cochain complexes is combinatorially defined by the Eulerian idempotents, the induced isomorphism of algebras on cohomologies
$$ \HH^{\ast}_{\oplus}(\mathcal{A},\mathcal{A}) \,\cong\,\HH^{\ast}_{\oplus}(E,E) $$
preserves the Hodge decomposition as well. It follows from Corollary \ref{KosDualCochain} and Corollary \ref{HodgeHHE} that there is an isomorphism of algebras preserving Hodge decomposition
$$ \HH^{\ast}(\U\mfa,\U\mfa)\,\cong\,\HH^{\ast}_{\oplus}(\mathcal{A},\mathcal{A}) \ .$$
By Lemma \ref{HodgeCup}, the cup product on $\HH^{\ast}(\U\mfa,\U\mfa)$ preserves the Hodge decomposition, i.e.,
\begin{equation} \la{cupdecomp} \HH^{(p),\ast}(\mfa) \cup \HH^{(q),\ast}(\mfa) \,\subseteq\, \,\HH^{(p+q),\ast}(\mfa)\ . \end{equation}
Finally, define the BV operator $\Delta\,:\,\HH^{\ast}(\U\mfa,\U\mfa) \rar \HH^{\ast-1}(\U\mfa,\U\mfa)$ by $\Delta:=\Psi^{-1} \circ B \circ \Psi$, where $B$ denotes the Connes differential. By \cite[Thm. 2.2]{BRZ}, $B$ decreases Lie-Hodge degree by $1$. It follows from Lemma \ref{hodgeduality} that $\Delta(\HH^{(p),\ast}(\U\mfa,\U\mfa))\,\subseteq\,\HH^{(p-1),\ast}(\U\mfa,\U\mfa))$. By \cite[Thm. 3.4.3]{G2},
$$ [a,b]_G\,=\, \Delta(a \cup b) -\Delta(a) \cup b -(-1)^{|a|} a \cup \Delta(b) \ .$$
It follows from \eqref{cupdecomp} that $[a,b]_G\,\in\,\HH^{(p+q-1),\ast}(\U\mfa,\U\mfa)$ for all $a\,\in\,\HH^{(p),\ast}(\U\mfa,\U\mfa)$, and for all $b\,\in\,\HH^{(q),\ast}(\U\mfa,\U\mfa)$.
\eproof

\bproof[Proof of Theorem \ref{genmain}]
Since the minimal $L_\infty$-model of $\mfa$ is finite dimensional, non-negatively graded and nilpotent, $\HH^{(p),\ast}(\mfa) \cup \HH^{(q),\ast}(\mfa)\,\subset \,\HH^{(p+q),\ast}(\mfa)$ by Proposition \ref{cupMinLinfty}. By \cite[Thm. 2.2]{BRZ},\\ $B(\HC^{(p)}_\ast(\mfa))\,\subset\,\HH^{(p-1)}_\ast(\mfa)$ for all $p \geq 1$. Since $\mfa$ is Koszul dual to a finite dimensional $C\,\in\,\cDGC_{k/k}$ that is equipped with a non-degenerate cyclic pairing, $\Psi^{-1}(\HH^{(p-1)}_\ast(\mfa))\,=\,\HH^{(p-1),\ast}(\mfa)$ by Lemma \ref{hodgeduality}. Thus,
\begin{eqnarray*} I\big[\Psi[\Psi^{-1}(B(\HC^{(p)}_\ast(\mfa))) \cup \Psi^{-1}(B(\HC^{(q)}_\ast(\mfa)))]\big] &\subset & I\big[\Psi[\HH^{(p-1),\ast}(\mfa) \cup \HH^{(q-1),\ast}(\mfa)]\big]\\
                                                                                                     &\subset & I\big[\Psi[\HH^{(p+q-2),\ast}(\mfa)]\big]\\
                                                                                                     & = & I[\HH^{(p+q-2)}_\ast(\mfa)]\\
                                                                                                     & \subset & \HC^{(p+q-2)}_\ast(\mfa)
                                                                                                     \end{eqnarray*}
The equality follows from Lemma \ref{hodgeduality} and the last inclusion is by \cite[Thm. 2.2]{BRZ}. The desired theorem now follows immediately from Proposition \ref{poissoncup}.
\eproof

Assume that $\mfa\,\in\,\DGL_k^+$ satisfies the conditions of Theorem \ref{genmain}. Further assume that a finite dimensional $C\,\in\,\cDGC_{k/k}$ that is Koszul dual to $\mfa$ carries a non-degenerate cyclic pairing. By Theorem \ref{liestronhom}, there is an action of $\rHC_\ast(\U\mfa)$ on $\rHH_\ast(\U\mfa)$, making the latter a graded Lie module over the former. The following result strengthening Theorem 3.4 of \cite{BRZ} holds under the above conditions.
\bthm \la{MainActiononHH}
For all $p,q \geqslant 1$, $\{\HC^{(p)}_\ast(\mfa),\HH^{(q)}_\ast(\mfa)\}\,\subseteq\,\HH^{(p+q-2)}_\ast(\mfa)$.
\ethm
\bproof
Note that
\begin{eqnarray*}
\Psi\big[\Psi^{-1}(B(\HC^{(p)}(\mfa)))\,,\,\Psi^{-1}(\HH^{(q)}(\mfa))\big]_G
& \subset & \Psi\big[\HH^{(p-1),\ast}(\mfa)\,,\,\HH^{(q),\ast}(\mfa) \big]_G\\
& \subset & \Psi\big(\HC^{(p+q-2),\ast}(\mfa)\big)\\
& = & \HH^{(p+q-2)}_\ast(\mfa)\ .
\end{eqnarray*}
The first inclusion above is by \cite[Thm. 2.2]{BRZ} and Lemma \ref{hodgeduality}. The second inclusion above is by Proposition \ref{cupMinLinfty}. The equality is by Lemma \ref{hodgeduality}. The desired theorem now follows immediately from Proposition \ref{actiononhh}.
\eproof
\subsection{Application to string topology} Let $X$ be a simply connected space of finite rational type. Recall (see \cite{FHT}) that associated to $X$ are a commutative cochain DG $\Q$-algebra $\mathcal{A}_X$, called the {\it Sullivan model} of $X$ and a connected chain DG Lie $\Q$-algebra $\mathfrak{a}_X$, called the {\it Quillen model} of $X$. The Quillen and Sullivan models determine the rational homotopy of $X$, and are Koszul dual to each other in the sense that there is a quasi-isomorphism of DG algebras
$$ \C^{\ast}(\mathfrak{a}_X; \Q) \stackrel{\sim}{\longrightarrow} \mathcal{A}_X\,,$$
where $ \C^* $ denotes the Chevalley-Eilenberg cochain complex with trivial coefficients.

Now, let $\mathcal{L}X = {\rm Map}(S^1, X) $ denote the free loop space over $X$. In the case when $X$ is a simply connected closed oriented manifold of dimension $d$, Chas and Sullivan \cite{ChS} constructed a product --- 
called the {\it loop product} --- on the (reduced) rational homology of $\LL X$:
$$ 
\bullet\,:\, \overline{\H}_{\ast}(\mathcal{L}X, \Q) \otimes \overline{\H}_{\ast}(\mathcal{L}X, \Q) \rar \overline{\H}_{\ast-d}(\mathcal{L}X, \Q)\ .$$
the free loop spaces $\mathcal{L}X$ is equipped with a natural circle action via rotation of loops. In addition to the usual homology, one may therefore consider the $S^1$-equivariant homology (resp., reduced $S^1$-equivariant homology) $\H^{S^1}(\mathcal{L}X, \Q)$ (resp., $\overline{\H}^{S^1}_\ast(\mathcal{L}X, \Q)$) of $\mathcal{L}X$. The two homology theories are related by the Gysin long exact sequence
\begin{equation}
\la{epgysin}
 \ldots\, \xrightarrow{D} \,\overline{\H}^{S^1}_{n-1}(\mathcal L X, \Q) \,\stackrel{\mathbb{M}}{\to}\,\overline{\H}_{n}(\mathcal L X, \Q) \,\stackrel{p_\ast}{\to}\, \overline{\H}^{S^1}_{n}(\mathcal L X, \Q) \,\xrightarrow{D}\, \overline{\H}^{S^1}_{n-2}(\mathcal L X, \Q)\,\to\, \ldots
\end{equation}
where $D$ stands for the Gysin map and $p\,:\,\LL X \times ES^1 \rar \LL X \times_{S^1} ES^1$ is the canonical projection. The {\it string bracket}  on $\rH^{S^1}(\mathcal{L}X, \Q)$ is the bilinear map induced by the loop product ({\it cf.} \cite{ChS}):
$$ \{\mbox{--},\mbox{--}\}\,:\,\overline{\H}^{S^1}(\mathcal{L}X, \Q) \otimes \overline{\H}^{S^1}(\mathcal{L}X, \Q) \rar \overline{\H}^{S^1}(\mathcal{L}X, \Q)\,,\qquad a \otimes b \mapsto \{a,b\}=(-1)^{|a|+d}p_\ast(\mathbb{M}(a) \bullet \mathbb{M}(b))\ . $$

The following theorem is a well-known result due to Goodwillie \cite{Go} and Jones \cite{J}
(see also \cite{JM}).
\bthm
\la{top1}
There are natural isomorphisms of graded vector spaces
\begin{equation*}
\la{ax}
\alpha_X:\, {\rHH}_{\ast}(\mathcal U\mathfrak{a}_X)  \xrightarrow{\sim} \overline{\H}_{\ast}(\mathcal L X, \Q)
\ ,\qquad
\beta_X:\, \rHC_{\ast}(\mathcal U\mathfrak{a}_X) \xrightarrow{\sim} {\rH}^{S^1}_{\ast}(\mathcal L X, \Q)
\end{equation*}
identifying the Connes periodicity sequence for $ \, \U\mathfrak{a}_X $ with the Gysin long exact sequence for the $S^1$-equivariant homology of $ \LL X $.
\ethm

On the other hand, the finite coverings of the circle $\varphi^n\,:\,S^1 \rar S^1 \,,\,\, e^{i\theta} \mapsto e^{ni\theta}$, give natural maps $\varphi^n_X\,:\,\LL X \rar \LL X$, one for each $ n \ge 0 $,
which induce Frobenius (power) operations on homology:
$$ \Phi^n_X\,:\,\overline{\H}_\ast(\LL X, \Q) \rar \overline{\H}_\ast(\LL X, \Q)\,, \qquad \tilde{\Phi}^n_X\,:\, \overline{\H}^{S^1}_\ast(\LL X, \Q) \rar  \overline{\H}^{S^1}_\ast(\LL X, \Q)\ .$$
By \cite[Theorem 4.1]{BRZ}, the isomorphisms $ \alpha_X $ and $ \beta_X $ of Theorem~\ref{top1}  restrict to isomorphisms of (graded) vector spaces
$$ {\HH}^{(p)}_\ast(\mathcal U\mathfrak{a}_X)  \xrightarrow{\sim} \overline{\H}^{(p)}_\ast(\LL X, \Q)\,, \qquad \HC^{(p)}_{\ast}(\mathfrak{a}_X) \xrightarrow{\sim} {\rH}^{S^1,(p-1)}_{\ast}(\mathcal L X, \Q)\,,$$
where the targets are common eigenspaces of the endomorphisms $ \Phi^n_X$ and $\tilde{\Phi}^n_X$ with eigenvalues $n^p$:
\begin{equation*}
{\rH}^{(p)}_{\ast}(\LL X, \Q)\, := \bigcap_{n \ge 0}\,\Ker(\Phi_X^n - n^p\,\id)\ ,\qquad
{\rH}^{S^1,\, (p)}_{\ast}(\LL X, \Q)\, :=
 \bigcap_{n \ge 0}\,\Ker(\tilde{\Phi}_X^n - n^p\,\id)
\end{equation*}
Thus, we have the Hodge-type decompositions
$$ \overline{\H}_\ast(\LL X, \Q)\,=\, \bigoplus_{p=0}^{\infty} \overline{\H}^{(p)}_\ast(\LL X, \Q)\,, \qquad  \overline{\H}^{S^1}_\ast(\LL X, \Q)\, =\, \bigoplus_{p=0}^{\infty}  \overline{\H}^{S^1,(p)}_\ast(\LL X, \Q)
$$
and the Gysin long exact sequence \eqref{epgysin} decomposes into a direct sum of exact sequences
\begin{equation} \la{hodgepgysin}
 \ldots\, \xrightarrow{D} \,\overline{\H}^{S^1,\, (p+1)}_{n-1}(\mathcal L X, \Q) \,\stackrel{\mathbb{M}}{\to}\,\overline{\H}^{(p)}_{n}(\mathcal L X, \Q) \,\stackrel{p_\ast}{\to}\, \overline{\H}^{S^1, (p)}_{n}(\mathcal L X, \Q) \,\xrightarrow{D}\, \overline{\H}^{S^1, \, (p+1)}_{n-2}(\mathcal L X, \Q)\,\to\, \ldots\ .
\end{equation}
By \cite[Theorem 4.2]{BRZ}, the string bracket gives $\overline{\H}^{S^1}_\ast(\LL X, \Q)$ the structure of a filtered Lie algebra with respect to the filtration
$$ F_p\overline{\H}^{S^1}_\ast(\LL X, \Q)\,:=\, \bigoplus_{q \leq p+1} \overline{\H}^{S^1,(q)}_\ast(\LL X, \Q) \ .$$
The following result, which strengthens \cite[Theorem 2]{FeT} in the rationally elliptic case.

\bthm \la{loopprodhodge}
Let $X$ be a simply connected closed oriented manifold of rationally elliptic type. Then, the loop product preserves the Hodge decomposition, i.e,
$$ \overline{\H}^{(p)}_\ast(\LL X, \Q) \bullet \overline{\H}^{(q)}_\ast(\LL X, \Q) \,\subseteq \, \overline{\H}^{(p+q)}_\ast(\LL X, \Q)\ .$$
\ethm
Let $\mathfrak{a}_X$ be a Quillen model of $X$. Recall the natural isomorphism $\alpha_X\,:\,\HH_\ast(\U\mathfrak{a}_X)\,\cong\,\H_\ast(\LL X, \Q)$ of  Theorem \ref{top1}. The following proposition is a consequence of \cite[Thm. D, Prop. 8]{FTV}
\bprop \la{loopprod}
The isomorphism $\alpha_X$ identifies the loop product on $\H_\ast(\LL X, \Q)$ with the product
$$ \bullet\,:\, \HH_\ast(\U\mathfrak{a}_X) \otimes \HH_\ast(\U\mathfrak{a}_X) \rar \HH_\ast(\U\mathfrak{a}_X)\,,\qquad a \otimes b \mapsto a \bullet b= \Psi(\Psi^{-1}(a) \cup \Psi^{-1}(b))\ .$$
\eprop
\bproof
Let $d:=\dim X$. Recall that the map $\alpha_X^{\ast}$ gives a natural isomorphism $\H^{\ast}(\LL X, \Q) \,\cong\, \H^{\ast}(\mfa_X;\U\mfa_X^{\vee})$, where $\U\mfa_X^{\vee}$ denotes the graded linear dual of $\mfa_X$, both sides of which we identify. By \cite[Theorem D]{FTV}, the isomorphism $[X] \cap \mbox{--}\,:\,\H^{d-\ast}(\mfa_X;\U\mfa_X^{\vee}) \rar \H_{\ast}(\mfa_X;\U\mfa_X^{\vee})$ transforms the coproduct on the right hand side to the dual of the loop product on $\H^{\ast}(\mfa_X;\U\mfa_X^{\vee})$. The graded linear dual of this isomorphism therefore, transforms the loop product on $\H_\ast(\LL X, \Q)$ to the product on $\H^{\ast}(\mfa_X;\U\mfa_X)$. Let $C$ denote the Lambrechts Stanley model \cite{LS} of $X$. By (the proof of) Lemma \ref{hodgeduality}, the isomorphism $[X] \cap \mbox{--}$ is induced by the map  $\phi \otimes \id_{\U\mfa_X^{\vee}}\,:\,E \otimes \U\mfa_X^{\vee} \rar C[-d] \otimes \U\mfa_X^{\vee}$, where $E:=C^{\ast}$ and $\phi\,:\,E \rar C[-d]$ is the isomorphism induced by the cyclic pairing on $C$. It follows that the map $\phi^\ast\,:\,E[d] \rar C$ coincides with $\phi[d]$, whence the graded linear dual of the map $[X] \cap \mbox{--}\,:\,\H^{d-\ast}(\mfa_X;\U\mfa_X^{\vee}) \rar \H_{\ast}(\mfa_X;\U\mfa_X^{\vee})$  coincides with the map $[X] \cap \mbox{--}\,:\,\H^{d-\ast}(\mfa_X;\U\mfa_X) \rar \H_{\ast}(\mfa_X;\U\mfa_X)$. By Lemma \ref{hodgeduality}, this in turn, is identified with the isomorphism $\Psi\,:\,\HH^{d-\ast}(\U\mfa_X,\U\mfa_X) \rar \HH_\ast(\U\mfa_X)$. Since the product on $\H^{\ast}(\mfa_X;\U\mfa_X)$ is identified with the cup product on $\HH^\ast(\U\mfa_X,\U\mfa_X)$ by Proposition \ref{HHCup}, the desired proposition follows.
\eproof
\bproof[Proof of Theorem \ref{loopprodhodge}]
Since the minimal $L_{\infty}$-model of $\mfa_X$ is finite-dimensional and nilpotent, Proposition \ref{cupMinLinfty} applies to $\mfa_X$. By Proposition \ref{cupMinLinfty} and Lemma \ref{hodgeduality},
$$ \HH^{(p)}_\ast(\mfa_X) \bullet   \HH^{(q)}_\ast(\mfa_X)\,\subseteq\, \HH^{(p+q)}_\ast(\mfa_X)\ .$$
The desired result is therefore, immediate from Proposition \ref{loopprod} and \cite[Theorem 4.1]{BRZ}, by which $\alpha_X$ identifies $\HH^{(p)}_\ast(\mfa_X)$ with $\H^{(p)}_\ast(\LL X, \Q)$ for all $p$.
\eproof
The following lemma, which is known to experts, is an immediate consequence of Theorem \ref{top1}, Proposition \ref{poissoncup} and Proposition \ref{loopprod}.
\blemma 
\la{StringPoiss}
The isomorphism $\beta_X\,:\,\rHC_{\ast}(\mathcal U\mathfrak{a}_X) \xrightarrow{\sim} {\rH}^{S^1}_{\ast}(\mathcal L X, \Q)$ identifies the string bracket on ${\rH}^{S^1}_{\ast}(\mathcal L X, \Q)$ with the derived Poisson bracket on $\rHC_{\ast}(\mathcal U\mathfrak{a}_X)$ induced by the Poincar\'{e} duality pairing on its Koszul dual.
\elemma
We are now in position to give a proof of our first theorem stated in the Introduction.

\bproof[Proof of Theorem~\ref{MainTheorem}]
Since the minimal $L_\infty$-model of $\mfa_X$ is finite dimensional and nilpotent, and since the Koszul dual of $\mfa_X$ equipped with the Poincar\'{e} duality pairing has a finite-dimensional model (namely, the graded linear dual of the Lambrechts-Stanley model), Theorem \ref{MainTheorem} follows from Lemma \ref{StringPoiss}, Theorem \ref{genmain} and \cite[Theorem 4.1]{BRZ}.
\eproof

By Theorem \ref{liestronhom}, the Poincar\'{e} duality pairing on the Koszul dual of $\mfa_X$ also induces an action of $\rHC_\ast(\U\mfa_X)$ on $\HH_\ast(\U\mfa_X)$, making the latter a graded Lie module over the former. By Theorem \ref{top1}, we have a graded Lie action of $\rH_\ast^{S^1}(\LL X, \Q)$ (with string topology bracket) on $\H_\ast(\LL X, \Q)$. The following result strengthens Theorem 4.3 $(ii)$ of \cite{BRZ} in the case when $X$ is rationally elliptic.

\bthm \la{hs1actononh}
Assume that $X$ is rationally elliptic. Then,
$$
\{\H^{S^1,(p)}_\ast(\LL X, \Q)\,,\,\H^{(q)}_\ast(\LL X, \Q)\} \,\subseteq\, \H^{(p+q-2)}_\ast(\LL X, \Q)\ . 
$$
\ethm
\bproof
Since the minimal $L_{\infty}$-model of $\mfa_X$ is finite-dimensional and nilpotent,  and since the Koszul dual of $\mfa_X$ equipped with the Poincar\'{e} duality pairing has a finite-dimensional model (namely, the graded linear dual of the Lambrechts-Stanley model), the desired result follows immediately from Theorem \ref{MainActiononHH} and \cite[Theorem 4.1]{BRZ}.
\eproof


\begin{thebibliography}{}
%
%
%
\bibitem{AH}
L. Avramov and S. Halperin,  \textit{Through the looking glass: a dictionary between
rational homotopy theory and local algebra}, Lecture Notes in Math. \textbf{1183} (1986), 1--27.
%
\bibitem{BK}
C. Barwick and D. Kan, {\it Relative categories: another model for the homotopy theory of homotopy categories}, Indag. Math. (N.S.) \textbf{23} (2012), 42--68.
%
\bibitem{BCER}
Yu.~Berest, X.~Chen, F.~Eshmatov and A.~Ramadoss, \textit{Noncommutative Poisson structures, derived representation schemes and Calabi-Yau algebras}, Contemp. Math. \textbf{583} (2012), 219--246.
%
\bibitem{BFPRW}
Yu.~Berest, G.~Felder, S.~Patotski, A.~C.~Ramadoss and T.~Willwacher, \textit{Representation homology, Lie algebra cohomology and the derived Harish-Chandra homomorphism}, J. Eur. Math. Soc. \textbf{19} (2017), no. 9, 2811--2893.
%
\bibitem{BFR}
Yu. Berest, G. Felder and A. Ramadoss, \textit{Derived representation schemes and noncommutative geometry}. Contemp. Math. \textbf{607} (2014), 113--162.
%
\bibitem{BKR}
Yu. Berest, G. Khachatryan and A. Ramadoss, \textit{Derived representation schemes and cyclic homology}, Adv. Math. \textbf{245} (2013), 625--689.
%
\bibitem{BR}
Yu. Berest and A. Ramadoss, \textit{Stable representation homology and Koszul duality}, J. Reine Angew. Math. \textbf{715} (2016), 143--187.
%
\bibitem{BRZ}
Yu.~Berest, A.~Ramadoss and Y. Zhang, \textit{Dual Hodge decompositions and derived Poisson brackets},
Selecta Math. \textbf{23} (2017),  2029--2070.
%
\bibitem{BL}
R.~Bocklandt and L.~Le Bruyn, {\it Necklace Lie algebras and noncommutative symplectic geometry}, Math. Z. \textbf{240} (2002), no. 1, 141-167.
%
\bibitem{BHM}
M. B\"{o}kstedt, W. C. Hsiang and I. Madsen, {\it The cyclotomic trace and algebraic K-theory of spaces}, Invent. Math. \textbf{111} (1993), no. 3, 465--539.
%
\bibitem{BFG91}
D.~Burghelea, Z.~Fiedorowicz and W.~Gajda, {\it Adams operations in Hochschild and cyclic homology of de Rham algebra and free loop spaces}, K-Theory \textbf{4} (1991), no. 3, 269--287.
%
\bibitem{BFG}
D.~Burghelea, Z.~Fiedorowicz and W.~Gajda, \textit{Power maps and epicyclic spaces}, J. Pure Appl. Alg. \textbf{96} (1994), no. 1, 1-14.
%
\bibitem{CR}
D.~Calaque and C.~Rossi, {\it Lectures on Duflo isomorphisms in Lie algebras and complex geometry}, EMS Series of Lectures in Mathematics. European Mathematical Society (EMS) Z\"{u}ruch (2011), viii+106 pp. ISBN: 978-3-03719-096-8.
%
%
\bibitem{ChS}
M. Chas and D. Sullivan, \textit{String topology}, arxiv preprint math.GT/9911159.
%
\bibitem{CEEY}
X.~Chen, A.~Eshmatov, F.~Eshmatov and S.~Yang, {\it The derived non-commutative Poisson bracket on Koszul Calabi-Yau algebras}, J. Noncommut. Geom. \textbf{11} (2017), no. 1, 111-160.
%
%
\bibitem{CB}
W. Crawley-Boevey, \textit{Poisson structures on moduli spaces of representations}, J. Algebra \textbf{325} (2011), 205-215.
%
\bibitem{deVV}
L. de Thanhoffer de V\"{o}lcsey and M. Van den Bergh, \textit{Calabi-Yau deformations and negative cyclic homology}, J. Noncommu. Geom. \textbf{12} (2018), no. 4, 1255--1291.
%
\bibitem{VD}
V. A. Dolgushev, {\it A proof of Tsygan's formality conjecture for an arbitrary smooth manifold.} Thesis (Ph.D.)–Massachusetts Institute of Technology. 2005. {\it ProQuest LLC}.
%
\bibitem{Dr}
V. G. Drinfeld, \textit{On quasitriangular quasi-Hopf algebras and on a group that is closely connected with $\mathrm{Gal}(\overline{\mathbb Q}/{\mathbb Q})$}, Leningrad Math. J.,
\textbf{2} (1991), 829--860.
%
\bibitem{Du}
M. Duflo, \textit{Caract\`{e}res des groupes et des alg\`{e}bres de Lie r\'{e}solubles},
Ann. Sci. Ecole Norm. Sup. \textbf{3} (1970), 23--74.
%
\bibitem{DHKS}
W.~Dwyer, P.~Hirschhorn, D.~Kan and J.~Smith, {\it  Homotopy Limit Functors on
Model Categories and Homotopical Categories}, Mathematical Surveys and Monographs \textbf{113}, AMS, Providence, RI, 2004.
%
\bibitem{FT}
B.~Feigin and B.~Tsygan, {\it Additive K-Theory and crystalline cohomology}, Funct. Anal. Appl. \textbf{19} (1985), no. 2, 124--132.
%
\bibitem{FHT}
Y.~Felix, S.~Halperin and J.-C.~Thomas,
\textit{Rational Homotopy Theory},
Graduate Texts in Mathematics \textbf{205}, Springer-Verlag, New York, 2001.
%
\bibitem{FeT}
Y.~Felix and J.-C.~Thomas, \textit{Rational BV algebra in string topology}, Bull. Soc. Math. France \textbf{136} (2008), no. 2, 311--327.
%
\bibitem{FTV}
Y.~Felix, J.-C.~Thomas and M. Vigu\'{e}-Poirrier, \textit{Rational string topology}, J. Eur. Math. Soc. \textbf{9} (2007), 123--156.
%
\bibitem{Ge}
E. Getzler, \textit{Lie theory for nilpotent $L_\infty$-algebras}, Ann. Math. \textbf{170}, (2009), no. 1, 271--301.
%
\bibitem{GK}
E. Getzler and M. Kapranov, \textit{Cyclic operads and cyclic homology}, Geometry, topology and physics, 167--201, Conf. Proc. Lecture Notes Geom. Topology, IV, Intl. Press, Cambridge MA 1995.
%
%
%
\bibitem{G}
V.~Ginzburg, {\it Noncommutative symplectic geometry, quiver varieties and operads}, Math. Res. Lett. \textbf{8} (2001), 377-400.
%
\bibitem{G2}
V.~Ginzburg, {\it Calabi-Yau algebras}, \texttt{arXiv:math.AG/0612139}.
%
\bibitem{Go}
Th. Goodwillie, \textit{Cyclic homology, derivations, and the free loopspace},
Topology \textbf{24} (1985), 187--215.
%
\bibitem{Her}
E.~Herscovich, \textit{Hochschild (co)homology and Koszul duality}, \texttt{arXiv:1405.2247}.
%
\bibitem{HMS}
D. Husemoller, J. Moore and J. Stasheff, \textit{Differential homological algebra and homogeneous spaces}, J. Pure Appl. Alg. \textbf{5} (1974), 113--185.
%
\bibitem{J}
J.D.S.Jones, \textit{Cyclic homology and equivariant homology},
Invent. Math. \textbf{87} (1987), 403--423.
%
\bibitem{JM}
J.D.S. Jones and J. McCleary, \textit{Hochschild homology, cyclic homology, and
the cobar construction}, London Math. Soc. Lecture Note Ser., \textbf{175},
Cambridge Univ. Press, Cambridge, 1992, pp. 53--65.
%
\bibitem{Ke}
B. Keller, \textit{Deriving DG categories}, Ann. Sci. \'{E}cole Norm. Sup. (4) \textbf{27} (1994), no. 1, 63--102.
%
\bibitem{Ke1}
B. Keller, \textit{A-infinity algebras, modules and functor categories}, Contemp. Math. \textbf{406} (2006), 67--93.
%
\bibitem{Ke2}
B. Keller, \textit{Derived invariance of higher structures on the Hochschild complex}, available at https://webusers.imj-prg.fr/~bernhard.keller/publ/dih.pdf
%
\bibitem{LS}
P.~Lambrechts and D.~Stanley, {\it Poincar\'{e} duality and commutative differential graded algebras}, Ann. Sci. \'{E}cole Norm. Sup. (4) \textbf{41}, No. 4 (2008), 497--511.
%
\bibitem{LSX}
H.-Y. Liao, M. Sti\'{e}non and P. Xu, {\it Formality theorem for differential graded manifolds}, C. R. Math. Acad. Sci. Paris \textbf{356} (2018), no. 1, 27--43.
%
\bibitem{L89}
J.-L. Loday, \textit{Op\'{e}rations sur l'homologie cyclique des alg\`{e}bres commutatives}, Invent. Math. \textbf{96} (1989), no. 1, 205--230.
%
\bibitem{L}
J.-L. Loday, {\it Cyclic homology}, Grundl. Math. Wiss. \textbf{301}, 2nd Ed., Springer-Verlag, Berlin, 1998.
%
\bibitem{LV}
J.-L. Loday and B. Vallette, \textit{Algebraic operads}. Grundlehren der Mathematischen Wissenschaften \textbf{346}. Springer, Heidelberg 2012. xxiv+634 pp.
%
\bibitem{Lu}
J. Lurie, \textit{Higher topos theory}, Annals of Mathematics Studies, \textbf{170}. Princeton University Press, Princeton, NJ, 2009.
%
\bibitem{Mo}
J.~C.~Moore, \textit{Differential homological algebra}, In: Actes du Congr\`{e}s International des Math\'{e}maticiens (Nice, 1970), Tome 1, Gauthier-Villars, Paris, 335--339 (1971).
%
\bibitem{N}
C.~Negron, \textit{The cup product on Hochschild cohomology via twisting cochains and applications to Koszul rings}, J. Pure Appl. Alg. \textbf{221} (2017), no. 5, 1112--1133.
%
\bibitem{PT}
M. Pevzner and C. Torossian, \textit{Isomorphisme de Duflo et cohomologie tangentielle},
J. Geom. Phys. \textbf{51} (2004), no. 4, 486--505.
%
\bibitem{Pir}
T.~Pirashvili, \textit{On the PROP corresponding to bialgebras}, Cah. Topol. G\'{e}om. Diff\'{e}r. Cat\'{e}g. \textbf{43} (2002), no. 3, 221--239.
%
\bibitem{Q1}
D. Quillen, \textit{Homotopical Algebra}, Lecture Notes in Math. \textbf{43}, Springer-Verlag, Berlin, 1967.
%
\bibitem{Q2}
D. Quillen, \textit{Rational homotopy theory}, Ann. of Math. (2) \textbf{90} (1969), 205--295.
%
\bibitem{R}
A.~C.~Ramadoss, \textit{The big Chern classes and the Chern character}. Internat. J. Math. \textbf{19} (2008), no. 6, 699--746.
%
\bibitem{RZ}
A.~C.~Ramadoss and Y.~Zhang, \textit{Cyclic pairings and derived Poisson structures}, New York J. Math. \textbf{25} (2019), 1--44.
%
\bibitem{TZ1}
T.~Tradler and M.~Zeinalian, \textit{Infinity structure of Poincar\'{e} duality spaces}, Appendix A by Dennis Sullivan, Alg. Geom. Topol. \textbf{7} (2007), 233--260.
%
\bibitem{TZ2}
T.~Tradler and M.~Zeinalian, \textit{Algebraic string operations}, K-Theory \textbf{38} (2007), no. 1, 59--82.
%
\bibitem{VdB}
M. Van den Bergh, \textit{Double Poisson algebras}, Trans. AMS \textbf{360} (2008), 5711-5769.
%
\bibitem{W}
C.~Weibel, \textit{An introduction to homological algebra}, Cambridge studies in advanced mathematics \textbf{38}, Cambridge University Press, 1994.
%
\end{thebibliography}
\end{document}